\documentclass[a4paper,11pt]{amsart}
\usepackage{amsmath,amssymb,amsthm,graphicx,amsfonts}
\usepackage{mathrsfs}
\usepackage{setspace}
\usepackage[all]{xy}

\usepackage[top=25mm,bottom=25mm,left=30mm,right=30mm]{geometry}
\usepackage{quiver}
\usepackage{tikz}
\usetikzlibrary{patterns}
\usepackage{pxpgfmark}
\usepackage{multirow}
\usetikzlibrary{intersections,arrows,calc,decorations.markings}
\usepackage[colorlinks, linkcolor=blue, anchorcolor=blue, citecolor=green]{hyperref}
\tikzset{->-/.style={decoration={  markings,  mark=at position #1 with
			{\arrow{>}}},postaction={decorate}}}
\tikzset{-<-/.style={decoration={  markings,  mark=at position #1 with
			{\arrow{<}}},postaction={decorate}}}

\theoremstyle{plain}          
\newtheorem{theorem}{Theorem}[section]
\newtheorem{proposition}[theorem]{Proposition}
\newtheorem{lemma}[theorem]{Lemma}
\newtheorem{corollary}[theorem]{Corollary}
\theoremstyle{definition}
\newtheorem{definition}[theorem]{Definition}
\newtheorem{remark}[theorem]{Remark}
\newtheorem{example}[theorem]{Example}

\def\Sf{\mathbf{S}}
\def\af{\mathbf{a}}
\def\lf{\mathbf{l}}
\def\Tf{\mathbf{T}}
\def\Pf{\mathbf{P}}
\def\Mf{\mathbf{M}}
\def\Mcal{\mathcal{M}}
\def\Ncal{\mathcal{N}}

\def\PAS{\mathbb{P}\mathbb{A}(\mathbf{S})}

\def\PTAS{\mathbb{P}\mathbb{A}^{\bowtie}(\mathbf{S})}
\def\<{\langle}
\def\>{\rangle}
\def\ba{\mathbf{a}}

\newcommand{\I}{\operatorname{(I)}}
\newcommand{\II}{\operatorname{(II)}}

\newcommand{\V}{\operatorname{(V)}}
\newcommand{\VI}{\operatorname{(VI)}}

\newcommand{\opname}[1]{\operatorname{#1}}

\newcommand{\Int}{\mathbf{Int}}
\newcommand{\Intv}{\underline{\mathbf{Int}}}
\newcommand{\dimv}{\underline{\mathbf{dim}}}

\newcommand{\Hom}{\mathrm{Hom}}

\newcommand{\nn}{node[black]{$\bullet$}}

\newcommand{\ie}{{\em i.e.}\ }

\title{Intersection vectors over skew-tilings}
\thanks{Partially supported by the National Natural Science Foundation of China (Grant No. 12171397, 12471037)}
\author[Deng]{Difan Deng}
\address{Difan Deng\\Department of Mathematics\\
	Southwest Jiaotong University\\
	610031 Chengdu \\
	P.R.China}
\email{difandeng@my.swjtu.edu.cn}

\author[Geng]{Shengfei Geng}
\address{Shengfei Geng\\Department of Mathematics\\
Sichuan University\\
610064 Chengdu\\
P.R.China}
\email{genshengfei@scu.edu.cn}

\author[Liu]{Pin Liu}
\address{Pin Liu\\Department of Mathematics\\
	Southwest Jiaotong University\\
	610031 Chengdu \\
	P.R.China}
\email{pinliu@swjtu.edu.cn}

\begin{document}
\begin{abstract}
We prove that under a mild condition, a multiset of tagged permissible arcs over a skew-tiling is uniquely determined by its intersection vector. As an application, it is proved that---up to isomorphism---different $\tau$-rigid modules over a skew-gentle algebra $A$ arising from a skew-triple $(Q,Sp,I)$ have different dimension vectors if and only if $(Q,I)$ has no  minimal oriented cycle of  even-length with full zero relations. This generalizes a recent work of Fu-Geng for gentle algebras.
\end{abstract}

\dedicatory{In memory of Idun Reiten.}
\maketitle	
\tableofcontents
\section{Introduction}\label{s:intro}

Cluster algebras were introduced by Fomin-Zelevinsky \cite{FZCA1}, with quiver mutation as the combinatorial aspect. Over time, these structures have appeared in various mathematical fields. A significant class arises from marked surfaces: Fomin-Shapiro-Thurston \cite{FST} constructed associated quivers via triangulations of such surfaces and showed that mutation of quivers is compatible with flip of triangulations. The categorification of cluster algebras \cite{BMRRT} further connects them to quiver representations, making marked surfaces a natural setting for cluster categories. 
In particular, by constructing skew-gentle algebras, Qiu-Zhou \cite{QZ17} showed that there is a bijection between tagged curves and string objects. Moreover the dimensions of $\operatorname{Ext}^1$ are interpreted as intersection numbers of tagged curves and Auslander-Reiten translation as tagged rotation.

Skew-gentle algebras, introduced by Gei{\ss}-de la Pe\~{n}a \cite{GCJ99}, form an important class of representation-tame finite-dimensional algebras. 
Gentle algebras, first defined in \cite{AS87}, have been extensively studied due to their close ties to derived categories and cluster algebras—particularly those arising from surface triangulations \cite{ABCP}. 
Baur-Sim\~{o}es \cite{BS} provided a geometric interpretation: gentle algebras are realized as tiling algebras, and their indecomposable modules correspond to permissible curves. Recently, He-Zhou-Zhu \cite{HZZ} extended this model to skew-gentle algebras by introducing skew-tiling algebras. They proved that every skew-gentle algebra arises in this way, and established a bijection between tagged permissible curves on skew-tiling and certain indecomposable modules over the corresponding skew-gentle algebra.


 In \cite{FG22}, Fu and the second author employed the geometric models from \cite{BS} to study the $\tau$-tilting theory of gentle algebras. They proved that a multiset of permissible arcs on a tiling is uniquely determined by its intersection vector under certain conditions. As an application, they established that---up to isomorphism---distinct $\tau$-rigid modules over a gentle algebra $A$ have distinct dimension vectors if and only if $A$ contains no even-length oriented cycle with full relations.


In this paper, we aim to extend this study to $\tau$-rigid modules over skew-gentle algebras using the geometric models developed in \cite{HZZ}. A {\it skew-tiling} $(\Sf,\Mf,\Pf,\Tf)$ is defined as a marked surface $(\Sf,\Mf,\Pf)$ with an admissible partial triangulation $\Tf$, where the surface is divided into tiles of types $\I\text{–}\VI$ (see Section \ref{s:skew-tiling} for details). Associated to such a skew-tiling with its admissible tagged partial triangulation $\Tf^{\bowtie}$, we construct the \emph{unfolded tiling} $(\Sf^*,\Mf^*,\Tf^*)$.
Let $\mathscr{R}(\Sf)$ denote the set of equivalence classes of  finite multisets of pairwise compatible tagged permissible arcs in $(\Sf,\Mf,\Pf,\Tf)$.
For a multiset $\Mcal\in \mathscr{R}(\Sf)$, we define an invariant  $S_{\Mcal}$ and associate a corresponding multiset $\Mcal^*\in\mathscr{R}(\Sf^*)$. The intersection vectors $\underline{\Int}_{\Tf^{\bowtie}}(\Mcal)$ and $\Intv_{\Tf^*}(\Mcal^*)$ encode geometric information about these arcs, and ultimately correspond to the dimension vectors of the associated $\tau$-rigid modules.

A key result establishes the equivalence between intersection vectors in the skew-tiling and its unfolded counterpart:
\begin{theorem}(Theorems~\ref {t:equivalence of T^* and T^bowtie} and \ref{t:tagged multiset determines})\label{t:equivalence of unfolded and skew}
    Let $(\Sf,\Mf,\Pf,\Tf)$ be a skew-tiling and $\Tf^{\bowtie}$ the tagged version of $\Tf$. 
Let $(\Sf^*,\Mf^*,\Tf^*)$ be the associated unfolded tiling of $(\Sf,\Mf,\Pf,\Tf)$. Then for any $\Mcal,\Ncal\in\mathscr{R}(\Sf)$, the following are equivalent:
\begin{itemize}
    \item [(1)] $\Intv_{\Tf^{\bowtie}}(\Mcal)=\Intv_{\Tf^{\bowtie}}(\Ncal)$; 
    \item [(2)]  $\Intv_{\Tf^*}(\Mcal^*)=\Intv_{\Tf^*}(\Ncal^*)$ and  $S_{\Mcal}=S_{\Ncal}.$
\end{itemize}
Moreover, if  $(\Sf,\Mf,\Pf,\Tf)$ contains no tile of type $\II$ or even gons of type $\V$, then $(1)$ is further equivalent to 
\begin{itemize}
    \item [(3)]   ${\Mcal}={\Ncal}.$
\end{itemize}
\end{theorem}

  

Combining Theorem \ref{t:equivalence of unfolded and skew} with the geometric model from \cite{HZZ}, we obtain a precise characterization of skew-gentle algebras whose $\tau$-rigid modules are uniquely determined by their dimension vectors:
\begin{theorem}(Theorem~\ref{t: ns-conditions of tau-rigid})\label{t:tau rigid over skew} 
   Let $(Q,Sp,I)$ be a skew-gentle triple,   
$A^{sg}=KQ^{sp}/\langle I^{sg}\rangle$
 the associated skew-gentle algebras. The following are equivalent:
    \begin{itemize}
        \item[(1)] Different $\tau$-rigid $A^{sg}$-modules have different dimension vectors;
        \item[(2)]  The determinant of the Cartan matrix of $A^{sg}$ is nonzero;
        \item[(3)]  $(Q,I)$ contains no minimal oriented cycle of  even-length with full zero relations;
         \item[(4)]  $(Q^{sp},I^{sg})$ contains no minimal oriented cycle of  even-length with full zero relations.
    \end{itemize}
 \end{theorem}
 
 The paper is organized as follows. In Section \ref{s:s2}, we provide necessary background on skew-gentle algebras and skew-tilings. Section \ref{s:s3} establishes a bijection from skew-tiling to its associated unfolded tiling (see Lemma \ref{l:skew-gentle to gentle}) and proves Theorem~\ref{t:equivalence of unfolded and skew}. 
   Section \ref{s:s4} is devoted to the proof of Theorem \ref{t:tau rigid over skew} using Theorem \ref{t:equivalence of T^* and T^bowtie} and the geometric interpretation of $\tau$-tilting theory for skew-gentle algebras. Finally, in Section~\ref{s:two examples}, we provide two examples.
\subsection*{Acknowledgements}
The authors thank Changjian Fu for his many helpful suggestions.

\section{Skew-tiling and   unfolded tiling}\label{s:s2} 
Throughout this paper, $K$ denotes an algebraically closed field. By a $K$-algebra, we mean a basic finite dimensional associative $K$-algebra. For an algebra $A$, we write $\mathrm{mod} A$ for the category of finitely generated right $A$-modules. Whenever we speak of different modules, it is always up to isomorphism.

Let $Q=(Q_0,Q_1,s,t)$ be a quiver, where $Q_0$ is the set of vertices, $Q_1$ is the set of arrows, and $s, t: Q_1\to Q_0$ are the source and target functions.  In addition, to each vertex $i\in Q_0$, we associate a trivial path  $\varepsilon_i$ of length $0$ with $s(\varepsilon_i)=i=t(\varepsilon_i)$. For two arrows $\alpha$ and $\beta$, we write their composition as $\alpha\beta$ if $t(\beta)=s(\alpha)$.
Let $\langle I\rangle$ be an ideal in the path algebra $KQ$, generated by a set of relations $I$. 

\subsection{Skew-gentle algebras and skewed-gentle algebras}\label{ss:Skew-gentle algebras and skewed-gentle algebras}
This subsection recalls the notion of skew-gentle algebras and  skewed-gentle algebras from \cite{G99,GCJ99,HZZ}. 

A pair $(Q,I)$ is called a \emph{gentle pair} if it satisfies:
	\begin{enumerate}
		\item[(G1)] Each vertex is the source of at most two arrows and the target of at most two arrows.
		\item[(G2)] The ideal $I$ is formed by zero relations of length 2.
		\item[(G3)] For any arrow $\alpha$, there is at most one arrow $\beta$ (resp. $\gamma$) such that $\beta\alpha\in I$ (resp. $\gamma\alpha\notin I$).
		\item[(G4)] For any arrow $\alpha$, there is at most one arrow $\beta$ (resp. $\gamma$) such that $\alpha\beta\in I$ (resp. $\alpha\gamma\notin I$).
	\end{enumerate}
	A finite dimensional basic algebra $A$ is \emph{gentle} if  $A\cong KQ/\langle I\rangle$ for a gentle pair  $(Q,I)$. 
    

    Now let $(Q,I)$ be a gentle pair and $Sp\subseteq Q_0$ a  subset of \emph{special vertices}, the remaining vertices are \emph{ordinary}. For a triple $(Q,Sp,I)$, we consider the pair $(Q^{sp},I^{sp})$,  where $Q_0^{sp}=Q_0$, $Q_1^{sp}=Q_1\cup\{\epsilon_i\mid i\in Sp \}$ with $\epsilon_i$ a loop at $i$ and $I^{sp}=I\cup\{\epsilon_i^2\mid i\in Sp \}$. We call $\epsilon_i$, where $i\in Sp$, a \emph{special loop} and the elements in $Q_1$ \emph{ordinary arrows}.  A triple $(Q,Sp,I)$ as above is called \emph {skew gentle} if the corresponding pair $(Q^{sp},I^{sp})$ is gentle. For simplicity, we denote by $A^{g}=KQ/\langle I\rangle$ and $A^{sp}=KQ^{sp}/\langle I^{sp}\rangle$ for a skew-gentle triple $(Q,Sp,I)$.
	
 \begin{definition}
      A finite dimensional basic algebra is called \emph{skew-gentle} if it is isomorphic to $K Q^{sp}/\langle I^{sg}\rangle$ for a skew-gentle triple $(Q,Sp,I)$, where $I^{sg}=I\cup \{\epsilon_i^2-\epsilon_i\mid i\in Sp\}$. 
 \end{definition} 
 For a skew-gentle triple $(Q,Sp,I)$, denote by $A^{sg}=K Q^{sp}/\langle I^{sg}\rangle$ the associated skew-gentle algebra.
By definition, the algebra $A^{sg}$ is obtained from $KQ^{sp}/\langle I^{sp}\rangle$ by replacing the nilpotency condition $\epsilon_i^2=0$ with the idempotency condition $\epsilon_i^2=\epsilon_i$ for all special loops $\epsilon_i$ with $i\in Sp$. Moreover,  it is easy to get that the set
\[
\{\varepsilon_j\mid j\in Q_0{\setminus}Sp\}\cup\{\varepsilon_i-\epsilon_i,\,\epsilon_i\mid i\in Sp\}
\]
forms a complete set of primitive orthogonal idempotents of $A^{sg}$.   


\begin{definition}
    A finite dimensional basic algebra is called \emph{skewed-gentle} if it is isomorphic to $KQ^{\bowtie}/I^{\bowtie}$ for a skew-gentle triple $(Q,Sp,I)$, where 
 $(Q^{\bowtie}, I^{\bowtie})$ as follows:
\begin{itemize}
    \item $Q_0^{\bowtie}=\{Q_0(i)\mid i\in Q_0\}$, where $Q_0(i)=\begin{cases}
        i & \text{ if } i\notin Sp\\
        i^{-}\cup i^{+}& \text{ if } i\in Sp
    \end{cases}$;
    \vspace{0.1cm}
    \item $Q^{\bowtie}_1=\big\{(a,\alpha,b)\mid \alpha\in Q_1,a\in Q_0(s(\alpha)), b\in Q_0(t(\alpha))\big\}$;
    \vspace{0.1cm}
    \item $I^{\bowtie}=\Big\{\sum\limits_{b\in Q_0(s(\alpha))}\lambda_b(b,\alpha,c)(a,\beta,b)\mid \alpha\beta\in I, a\in Q_0(s(\beta)),c\in Q_0(t(\alpha))\Big\},$ where $\lambda_b=-1$ if $b=i^{-}$ for some $i\in Q_0$ and $\lambda_b=1$ otherwise.
\end{itemize}
\end{definition}
For a skew-gentle triple $(Q,Sp,I)$, denote by $A^{{\bowtie}}=KQ^{\bowtie}/I^{\bowtie}$ the associated skewed-gentle algebra. 
According to \cite{GCJ99}, $A^{sg}$ is Morita equivalent to $A^{{\bowtie}}$.

\subsection{Marked surface}
As in \cite{FST}, a \emph{marked surface} is a triple $(\Sf,\Mf,\Pf)$, where $\Sf$ is a compact oriented surface with nonempty boundary $\partial \Sf$, $\Mf\subset \partial \Sf$ is a finite set of marked points on the boundary and $\Pf\subset \Sf\setminus\partial \Sf$ is a finite set of punctures in the interior of $\Sf$.

A connected component of $\partial \Sf$ is called a \emph{boundary component} of $\Sf$. A boundary component $B$ of $\Sf$ is \emph{unmarked} if $\Mf\cap B=\emptyset$. A \emph{boundary segment} is the closure of a component of $\partial \Sf\setminus \Mf$. 

A \emph{curve} on a marked surface $(\Sf,\Mf,\Pf)$ is a continuous map $\gamma:~[0,1]\longrightarrow \Sf$ satisfying:
    \begin{enumerate}
        \item both $\gamma(0)$ and $\gamma(1)$ are in $\Mf\cup \Pf$ and $\gamma(t)\in \Sf\setminus(\partial \Sf\cup \Pf)$ for $0<t<1$;
        \item $\gamma$ is neither null-homotopic nor homotopic to a boundary segment. 
    \end{enumerate}

The \emph{inverse} of $\gamma$ is defined as $\gamma^{-1}(t)=\gamma(1-t),$ for $t\in[0,1]$. 
For two curves $\gamma_1$ and $\gamma_2$, $\gamma_1 \sim\gamma_2$ means that $\gamma_1$ is homotopic to $\gamma_2$ relative to $\{0,1\}$ $(${\em i.e.} fixing the endpoints$)$. Define an equivalence relation $\cong$ on the set of curves in $\Sf$ by $\gamma_1\cong\gamma_2$ if $\gamma_1\sim\gamma_2$ or $\gamma_1^{-1}\sim\gamma_2$. 
Throughout, curves are considered up to homotopy relative to endpoints and up to inversion.

\subsection{Admissible partial triangulation}  
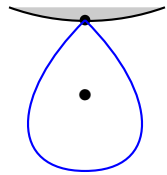
\begin{figure}[htpb]\centering
			\begin{tikzpicture}[xscale=1,yscale=1]
\draw[thick,fill=black!20] (2,0.17)arc (250:290:3);
\fill(3,0) circle(2pt);
\draw[thick,color=blue] (3,0) .. controls (2,-1) and (2,-2) .. (3,-2);
\draw[thick,color=blue] (3,0) .. controls (4,-1) and (4,-2) .. (3,-2);

\draw(3,-1)node{$\bullet$};
\end{tikzpicture}
		\caption{A once-punctured monogon}
		\label{fig:self-fold}
	\end{figure}

For convenience, through out this paper, by an {\it arc}, we always mean a curve who has no self-intersection in $\Sf\setminus(\Mf\cup \Pf)$.
Two arcs are called {\it compatible} if they have no intersections in $\Sf\setminus(\Mf\cup \Pf)$. 
A \emph{partial triangulation} of $\Sf$ is a collection of pairwise compatible arcs.
Recall from \cite{HZZ} an \emph{admissible partial triangulation} $\Tf$ is a partial triangulation in which every puncture is enclosed within a monogon of $\Tf$, as illustrated in Figure~\ref{fig:self-fold}.
The loop that cuts out a monogon and encloses exactly one puncture is simply called a {\it once-punctured} loop.
\subsection{Basic tiles}
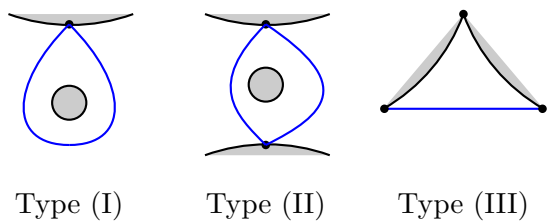
\begin{figure}[ht]
\begin{minipage}[t]{0.15\linewidth} 
\begin{tikzpicture}[xscale=0.8,yscale=0.8]

\draw[thick,fill=black!20] (2,0.17)arc (250:290:3);

\fill(3,0) circle(2pt);

\draw[thick,color=blue] (3,0) .. controls (2,-1) and (2,-2) .. (3,-2);
\draw[thick,color=blue] (3,0) .. controls (4,-1) and (4,-2) .. (3,-2);

\draw[thick,fill=black!20] (3,-1.3) circle(8pt);
\node at (3,-3){Type (I)};

\end{tikzpicture}
\end{minipage}
\begin{minipage}[t]{0.16\linewidth} 
\centering
		\begin{tikzpicture}[xscale=0.8,yscale=0.8]
\draw[thick,fill=black!20] (2,0.17)arc (250:290:3);

\fill(3,0) circle(2pt);

\draw[thick,fill=black!20] (2,-2.17)arc (110:70:3);
\fill(3,-2) circle(2pt);
\draw[thick,fill=black!20] (3,-1) circle(8pt);
\draw[thick,color=blue] (3,0) .. controls (2,-1) and (2.5,-1.5) .. (3,-2);
\draw[thick,color=blue] (3,0) .. controls (4.5,-1) and (4,-1.5) .. (3,-2);
\node at (3,-3){Type (II)};

\end{tikzpicture}
\end{minipage}
\begin{minipage}[t]{0.15\linewidth} 
\centering
		\begin{tikzpicture}[xscale=0.8,yscale=0.8]
\draw[thick,fill=black!20] (2,0.17)arc (200:240:3);
\draw[thick,fill=black!20] (2,0.17)arc (340:300:3);
\draw[blue,thick] (0.7,-1.4) to (3.3,-1.4);	
\fill(2,0.17) circle(2pt);
\fill(0.7,-1.4) circle(2pt);
\fill(3.3,-1.4) circle(2pt);
\node at (2,-3){Type (III)};

\end{tikzpicture}
\end{minipage}
\caption{Basic tiles of type (I)-(III)}\label{f:b1-3}
\end{figure}

\begin{figure}[ht]
   \begin{minipage}[t]{0.25\linewidth} 
\centering
		\begin{tikzpicture}[xscale=0.5,yscale=0.5]
					
	 \draw[thick,fill=black!40] (-1.3,-2.15)arc (120:60:1.5);
			
	\draw[blue,thick] (0,2) to (-1,2);
	\draw[blue,thick] (-1,2)\nn to (-2,1.5)\nn;
	\draw[blue,thick] (-2,1.5)\nn to (-2.7,0.5)\nn;
			
	\draw[blue,thick] (-1,-2)\nn to (-2,-1.5)\nn;
	\draw[blue,thick] (-2,-1.5)\nn to (-2.7,-0.5)\nn;
		
	\draw[blue,thick] (0,2)\nn to (1,1.5)\nn;
	\draw[blue,thick] (1,1.5)\nn to (1.7,0.5)\nn;
			
	\draw[blue,thick] (0,-2)\nn to (1,-1.5)\nn;
	\draw[blue,thick] (1,-1.5)\nn to (1.7,-0.5)\nn;
        \draw[blue](-2.7,-0.4)node[above]{$\vdots$}(1.7,-0.4)node[above]{$\vdots$};

	\node at (0,-3){ Type (IV)};
			\end{tikzpicture}
		
		\end{minipage}%
\begin{minipage}[t]{0.2\linewidth} 
\centering
		\begin{tikzpicture}[xscale=0.5,yscale=0.5]
		\draw[blue,thick] (0,2) to (-1,2);
		\draw[blue,thick] (-1,2)\nn to (-2,1.5)\nn;
		\draw[blue,thick] (-2,1.5)\nn to (-2.7,0.5)\nn;
	        \draw[blue,thick] (0,-2)\nn to (-1,-2)\nn;
		\draw[blue,thick] (-1,-2)\nn to (-2,-1.5)\nn;
		\draw[blue,thick] (-2,-1.5)\nn to (-2.7,-0.5)\nn;
		
		\draw[blue,thick] (0,2)\nn to (1,1.5)\nn;
		\draw[blue,thick] (1,1.5)\nn to (1.7,0.5)\nn;
			
		\draw[blue,thick] (0,-2)\nn to (1,-1.5)\nn;
		\draw[blue,thick] (1,-1.5)\nn to (1.7,-0.5)\nn;

		\draw[blue](-2.7,-0.4)node[above]{$\vdots$}(1.7,-0.4)node[above]{$\vdots$};

	\node at (0,-3){ Type (V)};
			\end{tikzpicture}	
	\end{minipage}	
  \begin{minipage}[t]{0.15\linewidth} 
\centering
		\begin{tikzpicture}[xscale=0.8,yscale=0.8]
\draw[thick,fill=black!20] (2,0.17)arc (250:290:3);

\fill(3,0) circle(2pt);

\draw[thick,color=blue] (3,0) .. controls (2,-1) and (2,-2) .. (3,-2);
\draw[thick,color=blue] (3,0) .. controls (4,-1) and (4,-2) .. (3,-2);

\draw(3,-1)node{$\bullet$};
\node at (3,-3){Type (VI)};
\end{tikzpicture}
\end{minipage}
    \caption{Basic tiles of type (IV)-(VI)} \label{f:b4-6}
\end{figure}
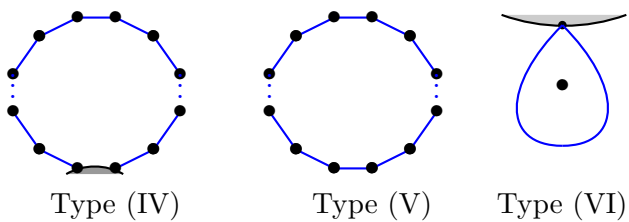

Let $\Tf$ be an admissible partial triangulation of a marked surface $(\Sf,\Mf,\Pf)$. Then $\Sf$ is divided by $\Tf$ into a collection of regions, also called \emph{tiles}. The following types of tiles are of particular importance (see Figures \ref{f:b1-3} and \ref{f:b4-6}):
\begin{enumerate}
	\item[(I)] monogons containing exactly one unmarked boundary component and no punctures in their interior;
	\item[(II)] digons containing exactly one unmarked boundary component and no punctures in their interior;
    \item [(III)]three-gons bounded by two boundaries and one arc in $\Tf$, and containing no unmarked boundary components or punctures in their interior;
    \item [(IV)] $m$-gons, whose edges are arcs in $\Tf$ and one boundary segment, and containing no unmarked boundary components or punctures in their interior; 
      \item[(V)] $m$-gons ($m\geq 3$), whose edges are arcs in $\Tf$ and containing no unmarked boundary components or punctures in their interior;
	\item[(VI)] once-punctured monogons, \ie monogons containing exactly one puncture and no unmarked boundary components in their interior.
\end{enumerate}

Note that the number of tiles of type $\VI$ is $|\Pf|$.
\subsection{Skew-tiling algebras}\label{s:skew-tiling}
Building on the geometric setup established above, we recall the definition of skew-tiling algebras, following \cite{HZZ}. 
\begin{definition}\label{def:st}
	Let $\Tf$ be an admissible partial triangulation of a marked surface $(\Sf,\Mf,\Pf)$. The \emph{skew-tiling algebra} associated to $\Tf$ is defined as $A_{\Tf}=K Q^\Tf/\langle R^\Tf\rangle$, where the quiver $Q^\Tf$ and the relation set $R^\Tf$ are constructed as follows:
	\begin{itemize}
		\item The vertices in $Q^\Tf_0$ correspond bijectively to the arcs in $\Tf$.
		\item An arrow from $i$ to $j$ exists if the corresponding arcs  $i$ and $j$ share an endpoint $p_{\alpha}\in \Mf$, and $j$ immediately follows $i$ in the anticlockwise order around $p_\alpha$. In particular, each loop in $Q^\Tf$ is at (the vertex indexed by) an arc in $\Tf$ whose endpoints coincide.    
		\item The relation set $R^\Tf$ includes:
		\begin{itemize}
			\item[(R1)] $\epsilon^2-\epsilon$ (resp. $\epsilon^2$) if $\epsilon$ is a loop such that the arc (corresponding to) $s(\epsilon)=t(\epsilon)$ does (resp. does not) cut out a region of type $\VI$; and
			\item[(R2)] $\alpha\beta$ if $p_{\beta}\neq p_{\alpha},$ 
            or the endpoints of the arcs (corresponding to) $s(\alpha)=t(\beta)$ coincide and the configuration matches one of the cases in Figure~\ref{fig:R2}.
			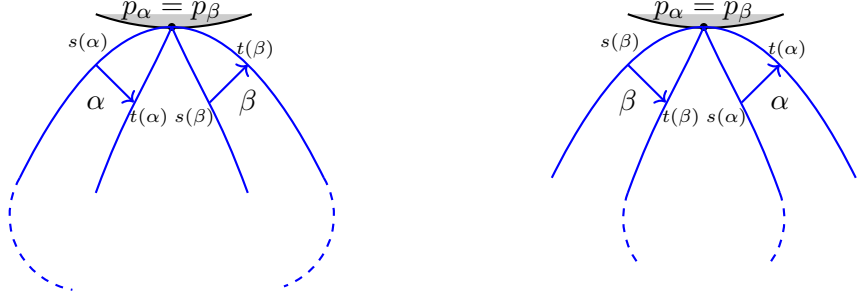
\begin{figure}[htbp]\centering
				\begin{tikzpicture}
\draw[thick,fill=black!20] (2,0.17)arc (250:290:3);
\node at (3,0.2){$p_\alpha=p_\beta$};
\fill(3,0) circle(1.5pt);
\draw[thick,color=blue] (3,0) parabola (1,-2);
\draw[thick,color=blue,dashed] (1,-2) arc (150:260:1);
\draw[thick,color=blue,->] (2,-0.5)--(2.5,-1);
\node at (2.7,-1.2){\tiny{$t(\alpha)$}};
\node at (1.9,-0.2){\tiny{$s(\alpha)$}};
\node at (2,-1) {$\alpha$};
\draw[thick,color=blue] (3,0) .. controls (2.6,-0.9) and (2.4,-1.1) .. (2,-2.2);
\draw[thick,color=blue] (3,0) .. controls (3.4,-0.9) and (3.6,-1.1) .. (4,-2.2);
\draw[thick,color=blue] (3,0) parabola (5,-2);
\draw[thick,color=blue,dashed] (5,-2) arc (30:-70:1);
\draw[thick,color=blue,->] (3.5,-1)--(4,-0.5);
\node at (4,-1) {$\beta$};
\node at (3.3,-1.2){\tiny{$s(\beta)$}};
\node at (4.1,-0.3){\tiny{$t(\beta)$}};

\draw[thick,fill=black!20] (9,0.17)arc (250:290:3);
\node at (10,0.2){$p_\alpha=p_\beta$};
\fill(10,0) circle(1.5pt);
\draw[thick,color=blue] (10,0) parabola (8,-2);
\draw[thick,color=blue,dashed] (9,-2.2) arc (160:220:1);
\draw[thick,color=blue,->] (9,-0.5)--(9.5,-1);
\node at (9,-1) {$\beta$};
\node at (9.7,-1.2){\tiny{$t(\beta)$}};
\node at (8.9,-0.2){\tiny{$s(\beta)$}};
\draw[thick,color=blue] (10,0) .. controls (9.6,-0.9) and (9.4,-1.1) .. (9,-2.2);
\draw[thick,color=blue] (10,0) .. controls (10.4,-0.9) and (10.6,-1.1) .. (11,-2.2);
\draw[thick,color=blue] (10,0) parabola (12,-2);
\draw[thick,color=blue,dashed] (11,-2.2) arc (20:-40:1);
\draw[thick,color=blue,->] (10.5,-1)--(11,-0.5);
\node at (11,-1) {$\alpha$};
\node at (10.3,-1.2){\tiny{$s(\alpha)$}};
\node at (11.1,-0.3){\tiny{$t(\alpha)$}};
\end{tikzpicture}
\caption{Relations in (R2), the case $p_\beta=p_\alpha$ and $s(\alpha)=t(\beta)$ is a loop}
				\label{fig:R2}
			\end{figure}
		\end{itemize}
	\end{itemize}
\end{definition}
Note that when $\Pf$ is empty, any partial triangulation is admissible. In this case, skew-tiling algebras coincide with the tiling algebras introduced in \cite{BS}, and we have the following results:
\begin{theorem}\cite{BS}\label{t:tiling algebra}
 Let  $(\Sf,\Mf,\Pf)$ be a  marked surface with an admissible partial triangulation $\Tf$ and $\Pf=\emptyset$. Then the algebra $A_{\Tf}$ is gentle. Conversely, for any gentle algebra $A$, there exists an admissible partial triangulation $\Tf$ of a marked surface $(\Sf,\Mf,\Pf)$ with $\Pf =\emptyset$ such that $A_{\Tf}\cong A$, and $\Tf$ divides $(\Sf,\Mf,\Pf)$ into tiles of types $\I-\V$.
 \end{theorem}
\begin{theorem}\cite{HZZ}\label{t:skew-tiling algebra}
    Let $(\Sf,\Mf,\Pf)$ be a marked surface with an admissible partial triangulation $\Tf$. Then the skew-tiling algebra $A_{\Tf}$ is skew-gentle. Conversely, for any skew-gentle algebra $A$, there exists a marked surface $(\Sf,\Mf,\Pf)$ with an admissible partial triangulation $\Tf$ such that $A_{\Tf}\cong A$ and $\Tf$ divides $\Sf$ into tiles of types $\I-\VI$.	
\end{theorem}

 In the following, for an admissible partial triangulation $\Tf$ of $(\Sf,\Mf,\Pf)$, if 
  $\Tf$  divides $(\Sf,\Mf,\Pf)$ into tiles of types $\I-\VI$, we simply call $\Tf$  a {\it skew-tiling} of $(\Sf,\Mf,\Pf)$, or refer to $(\Sf,\Mf,\Pf,\Tf)$ as a {\it skew-tiling}. In particular, if $\Pf=\emptyset$, we write  
$(\Sf,\Mf,\Tf)$ and call it a {\it tiling}.

\subsection{Tagged  arcs}
  Let $(\Sf,\Mf,\Pf)$ be a  marked surface. A {\it tagged  arc} is an  arc in which each end is assigned one of two tags, plain or notched, subject to the following conditions:
\begin{itemize}
    \item[$(\mathbf{T1})$] It does not cut out a once-punctured monogon;
    \item [$(\mathbf{T2})$] Any end incident to the boundary marked point is tagged plain;
    \item  [$(\mathbf{T3})$]  Both ends of a loop are tagged in the same way.
\end{itemize}
For a tagged   arc $\gamma$, we define the following types:
\begin{itemize}
    \item  $\gamma$ is a \emph{plain arc} if both of its ends are tagged plain;
   \item  $\gamma$ is a \emph{1-notched arc} if one end is tagged plain and the other is tagged notched;
   \item $\gamma$ is a \emph{2-notched arc} if both of its ends are tagged notched.
\end{itemize}
 In the figures, the plain tags are typically omitted, while the notched tags are represented by the symbol $\bowtie$. The following figures illustrate these tagging conventions:
\begin{figure}[h]
\begin{minipage}[t]{0.3\linewidth} 
\begin{tikzpicture}[xscale=0.6,yscale=0.6]
\draw[blue,thick] (-4,0) to (-2,0)\nn;
\draw[blue,thick](2,0) to (4,0)\nn;
\node at (-5,0){plain};
\node at (0.5,0){notched};
\node at (3.5,0)[rotate=90]{$\bowtie$};
\end{tikzpicture}
\end{minipage}
\end{figure}

\begin{definition}
    Let $\alpha$ and $\beta$ be two tagged  arcs in $(\Sf,\Mf,\Pf)$ whose  untagged version coincide. If exactly one of $\alpha$ or $\beta$ is a 1-notched arc, then the pair $\{\alpha,\beta\}$ is called {\it a pair of conjugate arcs}. 
\end{definition}
  As a consequence of condition $(\mathbf{T3})$, if $\{\alpha,\beta\}$ is a pair of conjugate arcs, then neither $\alpha$ nor $\beta$ is a loop. 

\subsection{QZ-intersection number}
 Let $(\Sf,\Mf,\Pf)$ be a marked surface. The intersection number for tagged arcs was defined in \cite{QZ17}, following the approach introduced in \cite{FST}.
 
\begin{definition}{\cite[Definition 3.3]{QZ17}}\label{def:Int}
Let $\gamma, \delta$ be two tagged arcs. Their \emph{intersection number} $\Int(\gamma|\delta)$ is defined as $$\Int(\gamma|\delta)=\Int^A(\gamma|\delta)+\Int^{C}(\gamma|\delta)+\Int^D(\gamma|\delta),$$ where
\begin{itemize}
 \item[(1)] $\Int^A(\gamma|\delta)$ counts the number of intersections of $\gamma$ and $\delta$ in $\Sf\setminus(\Mf\cup \Pf)$;
 \item[(2)] $\Int^{C}(\gamma|\delta)=0$ unless $\gamma$ and $\delta$ form a pair of conjugate arcs, in which case $\Int^{C}(\gamma|\delta)=-1$;
 \item[(3)] $\Int^D(\gamma|\delta)$ counts the number of pairs of an end of $\gamma$ and an end of $\delta$ such that they are incident to a common puncture and their tags are different.
\end{itemize}  
\end{definition}

The intersection number in Definition \ref{def:Int} is symmetric but different from the ``intersection number'' $(\gamma | \delta)$ defined in \cite[Definition 8.4]{FST} (see \cite{Y24} for a comparison).


 \subsection{Admissible partial tagged triangulation}
 Two tagged arcs $\gamma$ and $\delta$ are \emph{compatible} if $\Int(\gamma|\delta)=0$. In particular, for any conjugate pair $\{\gamma,\gamma'\}$, $\gamma$ and $\gamma'$ are compatible.
 A \emph{partial tagged triangulation} of $\Sf$ is a collection of pairwise compatible tagged arcs. An \emph{admissible partial tagged triangulation}  is a partial tagged triangulation in which every puncture is connected to a boundary marked point by a pair of conjugate arcs and no other tagged arcs connect to this puncture.


Let $\Tf$ be an admissible partial triangulation on $(\Sf,\Mf,\Pf)$. For each $P\in\Pf$, we denote the unique loop in $\Tf$ that encloses $P$ by $\mathbf{l}_P$, denote by $b_P$ the base point of $\lf_P$. So there is  a unique pair    of conjugate arcs which can be  enclosed by $\lf_P$, denote this pair by $\{\mathbf{a}_P^-,\mathbf{a}_P^{+}\}$, where $\ba_P^-$ is the arc tagged plain at the end incident to $P$ while $\ba_P^{+}$ is the arc tagged notched at the end incident to $P$ ({\em cf.} Figure~\ref{f: Notation of lf_P}).

\begin{figure}[h]
	\begin{tikzpicture}[xscale=0.8,yscale=0.8]
\draw[thick,fill=black!20] (4,-3)arc (70:110:3);

\fill(3,-2.8) circle(2pt);
\draw[thick,color=blue] (3,-2.8) .. controls (1,-1) and (1,2) .. (3,2);
\draw[thick,color=blue] (3,-2.8) .. controls (5,-1) and (5,2) .. (3,2);

\draw(3,0)node{$\bullet$};
\draw [thick,blue](3,0)node[below]{$P$}(3,-2.8)node[below] {$b_P$}(3,2)node[below]{$\lf_P$};
\end{tikzpicture}
 \begin{tikzpicture}[xscale=0.8,yscale=0.8]
\draw[thick,fill=black!20] (4,-3)arc (70:110:3);

\draw[thick,color=blue,dashed] (3,-2.8) .. controls (1,-1) and (1,2) .. (3,2);
\draw[thick,color=blue,dashed] (3,-2.8) .. controls (5,-1) and (5,2) .. (3,2);

\fill(3,-2.8) circle(2pt);
\draw[thick,color=blue] (3,-2.8) .. controls (2.5,-2) and (2.5,-1) .. (3,0);
\draw[thick,color=blue] (3,-2.8) .. controls (3.5,-2) and (3.5,-1) .. (3,0);

\draw(3,0)node{$\bullet$};
\draw [thick,blue](3,0)node[above]{$P$}(3,2)node[below]{$\lf_P$};
 \node at (3.18,-0.4)[rotate=20]{$\bowtie$};
 \draw [thick,blue](2.7,-1)node[left]{$\ba_P^-$}(3.3,-1)node[right]{$\ba_P^{+}$}(3,-2.8)node[below] {$b_P$};
\end{tikzpicture}
\caption{ $\lf_P, b_P, \{\ba_P^{-},\ba_P^{+}\}$}
\label{f: Notation of lf_P}
\end{figure}

Let  $\Tf^{\bowtie}$ be the set of tagged arcs obtained from $\Tf$ by replacing each once-punctured loop with a corresponding pair of conjugate tagged arcs which can be enclosed by this loop, \ie $$\Tf^{\bowtie}=(\Tf\setminus\{\lf_P|P\in \Pf\})\cup \{\af_P^-,\af_P^+|P\in\Pf\}.$$
({\em cf.} Figure \ref{f:T to T^x}).
It is clear that 
$\Tf^{\bowtie}$ is an admissible partial tagged triangulation
on $(\Sf,\Mf,\Pf)$.
Then we call $\Tf^{\bowtie}$  the {\it tagged version of $\Tf$}. 
\begin{figure}[h]
	\begin{tikzpicture}[xscale=0.8,yscale=0.8]
\draw[thick,fill=black!20] (4,-3)arc (70:110:3);

\fill(3,-2.8) circle(2pt);
\draw[thick,color=blue] (3,-2.8) .. controls (1,-1) and (1,2) .. (3,2);
\draw[thick,color=blue] (3,-2.8) .. controls (5,-1) and (5,2) .. (3,2);

\draw(3,0)node{$\bullet$};
\draw [thick,blue](3,0)node[below]{$P$}(3,-2.8)node[below] {$b_P$}(3,2)node[below]{$\lf_P$};
\draw [thick,red,->](5,-1) to (6,-1);
\end{tikzpicture}
 \begin{tikzpicture}[xscale=0.8,yscale=0.8]
\draw[thick,fill=black!20] (4,-3)arc (70:110:3);

\fill(3,-2.8) circle(2pt);
\draw[thick,color=blue] (3,-2.8) .. controls (2.5,-2) and (2.5,-1) .. (3,0);
\draw[thick,color=blue] (3,-2.8) .. controls (3.5,-2) and (3.5,-1) .. (3,0);

\draw(3,0)node{$\bullet$};
\draw [thick,blue](3,0)node[above]{$P$};
 \node at (3.18,-0.4)[rotate=20]{$\bowtie$};
 \draw [thick,blue](2.7,-1.5)node[left]{$\ba_P^-$}(3.3,-1.5)node[right]{$\ba_P^{+}$}(3,-2.8)node[below] {$b_P$};
\end{tikzpicture}
\caption{From  $\Tf$ to $\Tf^{\bowtie}$}
\label{f:T to T^x}
\end{figure}
 
\subsection{Unfolded tiling}\label{ss:unfolded tiling}
Let $(\Sf,\Mf,\Pf)$ be a marked surface and $\Tf$ a skew-tiling with tagged version $\Tf^{\bowtie}$. 
Define a new marked surface $(\Sf^*,\Mf^*)$ by replacing each puncture $P$ with a boundary component $B_P$ containing a single marked point (still labeled $P$). Then $(\Sf^*,\Mf^*)$ has no punctures. 
Let $\Tf^*$ be the collection of arcs obtained from $\Tf^{\bowtie}$ by replacing $\af_P^{+}$ for each  pair of conjugate arcs $\{\af_P^-,\af_P^{+}\}$ with the  loop $\lf_P$ that closely wraps around $\alpha_P^{-}$ and $B_P$ ({\em cf.} Figure \ref{fig:T^bowtie to T*}), \ie 
$$\Tf^*=(\Tf^{\bowtie}\setminus\{\af_P^+|P\in\Pf \})\cup\{\lf_P\in \Tf|P\in \Pf\}.$$
Equivalently, 
$\Tf^{*}$ is obtained from $\Tf$ by adding, for each puncture $P$, a radius-like arc inside the once-punctured loop which encloses $P$, connecting $P$ to $b_P$, \ie $$\Tf^*=\Tf\cup\{\af_P^{-}\in \Tf^{\bowtie}|P\in \Pf\}$$
 ({\em cf.} Figures~\ref{f:T to T^x} and \ref{fig:T^bowtie to T*}). Clearly $\Tf^*$ is a partial triangulation of $(\Sf^*,\Mf^*)$, yielding a tiling $(\Sf^*,\Mf^*,\Tf^*)$.


\begin{figure}[ht]
\begin{tikzpicture}[xscale=0.8,yscale=0.8]
        \draw[thick,fill=black!20] (1,-2.15)arc (70:110:3);
\draw[thick,color=blue] (0,-2)\nn .. controls (-0.5,-1.5) and (-0.5,-0.5) .. (0,0)\nn;	
\draw[thick,color=blue] (0,-2) .. controls (0.5,-1.5) and (0.5,-0.5) .. (0,0);	
\draw[blue](0.2,-.25)node[rotate=20]{$\bowtie$};	
\draw[thick,blue](-0.3,-1)node[left]{$\af_P^{-}$};
\draw[thick,blue](0,0)node[above]{P}(0.3,-1)node[right]{$\af_P^{+}$}(0,-2)node[below]{$b_P$}(1,0)(0,0.8);
		\end{tikzpicture}
  \begin{tikzpicture}[xscale=0.8,yscale=0.8]
  
\draw[white,thick,=>] (0,2) to (2,2);
\draw[red,thick,->] (0,0) to (2,0);
\draw[white,thick,->] (0,-2) to (2,-2);
	\end{tikzpicture}
\begin{tikzpicture}[xscale=0.8,yscale=0.8]
 \draw[thick,fill=black!20] (1,-2.15)arc (70:110:3);
\draw[thick,fill=black!20] (0,0.45) circle(12pt);
\draw (0,-2)\nn  ;	
\draw[blue,thick,>=stealth] (0,-2)\nn to (0,0)\nn;
\draw[thick,color=blue] (0,-2) .. controls (-2,-1.5) and (-2,2) .. (0,2);	
\draw[thick,color=blue] (0,-2) .. controls (2,-1.5) and (2,2) .. (0,2);	
\draw[thick,blue](0,-1)node[left]{$\af_{P}^{-}$} (0,2)node[above]{$\lf_P$}(0,0)(0,-2)node[below] {$b_P$}(0,0.8)node[above]{$B_P$};
\draw[thick,blue](0,-0.5)node[right]{$P$};
\end{tikzpicture}
\caption{From  $\Tf^{\bowtie}$ to $\Tf^*$}
\label{fig:T^bowtie to T*}
\end{figure}
According to \cite{HZZ}, replacing each puncture in $(\Sf,\Mf,\Pf, \Tf)$ with an unmarked boundary component ({\em cf.} Figure~\ref{fig:T to T^{circ}}) gives another tiling $(\Sf_{0},\Mf_{0},\Tf_{0})$. We call $(\Sf_{0},\Mf_{0},\Tf_{0})$ the associated \emph{original tiling} and $(\Sf^*,\Mf^*,\Tf^*)$ the associated \emph{unfolded tiling}.
\begin{figure}[h]
	\begin{tikzpicture}[xscale=0.8,yscale=0.8]
\draw[thick,fill=black!20] (4,-3)arc (70:110:3);

\fill(3,-2.8) circle(2pt);
\draw[thick,color=blue] (3,-2.8) .. controls (2,1) and (2,2) .. (3,2);
\draw[thick,color=blue] (3,-2.8) .. controls (4,1) and (4,2) .. (3,2);

\draw(3,0)node{$\bullet$};
\draw [thick,blue](3,0)node[below]{$P$};
\end{tikzpicture}
 \begin{tikzpicture}[xscale=0.8,yscale=0.8]
  
\draw[white,thick,=>] (0,2) to (2,2);
\draw[red,thick,->] (0,0) to (2,0);
\draw[white,thick,->] (0,-2) to (2,-2);
	\end{tikzpicture}
\begin{tikzpicture}[xscale=0.8,yscale=0.8]
\draw[thick,fill=black!20] (4,-3)arc (70:110:3);
\draw[thick,fill=black!20] (3,0) circle(8pt);
\fill(3,-2.8) circle(2pt);
\draw[thick,color=blue] (3,-2.8) .. controls (2,1) and (2,2) .. (3,2);
\draw[thick,color=blue] (3,-2.8) .. controls (4,1) and (4,2) .. (3,2);


\end{tikzpicture}
\caption{From  $\Tf$ to $\Tf_{0}$}
\label{fig:T to T^{circ}}
\end{figure}
\subsection{Quiver with relations for skew-tiling and unfolded tiling}\label{ss:construction of skewed and unfolded}
Let  $(\Sf,\Mf,\Pf)$ be a marked surface with a skew-tiling $\Tf$.
Let $\Tf^{sp}$ be the subset of $\Tf$ consisting of once-punctured loops (\ie $\Tf^{sp}=\{\lf_P\in\Tf|P\in\Pf\}$), and $Sp$ the subset of $Q_0^{\Tf}$ corresponding to arcs in $\Tf^{sp}$.

Define a quiver $Q$ by setting $Q_0=Q^{\Tf}_0$ and $Q_1=Q_1^{\Tf}\setminus\{\epsilon\in Q_1^{\Tf}| s(\epsilon)=t(\epsilon)\in Sp\}$.
Let $$I=R^{\Tf}\setminus\{\epsilon^2-\epsilon|\epsilon\in Q_1^{\Tf}, s(\epsilon)=t(\epsilon)\in Sp\}.$$ 
Then $(Q,Sp,I)$ is a skew-gentle triple. Moreover, $$Q^{\Tf}=Q^{sp}, I^{sg}=R^{\Tf}, I^{sp}=I\cup \{\epsilon^2|\epsilon\in Q_1^{\Tf}, s(\epsilon)=t(\epsilon)\in Sp\},$$
and it is clear that  $$A_{\Tf}\cong KQ^{sp}/\langle I^{sg}\rangle,\ \ \ A_{\Tf_0}\cong KQ^{sp}/\langle I^{sp}\rangle.$$



\begin{figure}[ht]
\begin{tikzpicture}[xscale=0.8,yscale=0.8]
 \draw[thick,fill=black!20] (1,-2.15)arc (70:110:3);
\draw[thick,fill=black!20] (0,0.45) circle(12pt);
\draw (0,-2)\nn  ;	
\draw[blue,thick,>=stealth] (0,-2)\nn to (0,0)\nn;
\draw[thick,color=blue] (0,-2) .. controls (-2,-1.5) and (-2,2) .. (0,2);	
\draw[thick,color=blue] (0,-2) .. controls (2,-1.5) and (2,2) .. (0,2);	
\draw[thick,color=red, <-] (-0.8,-1.4) to (-0.1,-1.4);
\draw[thick,color=red, <-] (0.1,-1.4) to (0.8,-1.4);
\draw[thick,red](-0.4,-1.4)node[above]{$\rho_{i^*}$}(0.5,-1.4)node[above]{$\rho_{i}$};
\draw[thick,blue](.1,-0.5)node[left]{$i^*$} (0,2)node[above]{$i$};
\draw[thick,blue](0,0)(0,-2)node[below] {$b_P$}(0,0.8)node[above]{$B_P$};
\end{tikzpicture}
\caption{$\rho_i\rho_{i^*}=0$}
\label{fig: rho_irho_i*=0}
\end{figure}
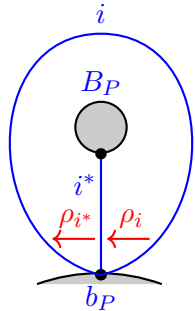

The next result provides a presentation of the unfolded tiling algebra in terms of quivers and relations.
\begin{proposition}\label{p:relation between tiling algebras}
The unfolded tiling algebra $A_{\Tf^*}$ is isomorphic to the path algebra $A^*$ of the bounded quiver $(Q^{*}, I^{*})$, where
\begin{itemize}
    \item $Q^*_0=Q_0\cup \{i^{*}\mid i\in Sp\}$;
     \vspace{0.1cm}
    \item $Q^*_1=Q_1\cup\{\rho_i:i\rightarrow i^*, \rho_{i^*}:i^*\rightarrow i\mid i\in Sp\}$;
     \vspace{0.1cm}
    \item $I^*=I\cup \{\rho_{i}\rho_{i^*}\mid i\in Sp\}$.
\end{itemize} 
\end{proposition}
\begin{proof}
 Consider the unfolded tiling $(\Sf^*,\Mf^*,\Tf^*)$. 
 Note that
 $$\Tf^*=\Tf\cup\{\af_P^{-}\in \Tf^{\bowtie}|P\in \Pf\}\text{  and  } \Tf^{sp}=\{\lf_P\in\Tf|P\in\Pf\}.$$
  For $i\in Sp$, let  $\lf_P\in \Tf^{sp}$ be the arc that indexes $i$ and $\af_P^{-}\in \Tf^*$ be the arc that indexes $i^*$. So $Q_0^{\Tf^*}=Q_0\cup\{i^*\mid i\in Sp\}$. Since $\af_P^{-}$ is closely wrapped by $\lf_P$, thus $$Q_1^{\Tf^*}=Q_1^{\Tf}\cup\{\rho_i,\rho_{i^*}\mid \rho_i:i\to i^*,\rho_{i^*}:i^*\to i,i\in Sp\}$$ and   $$R^{\Tf^*}=(R^\Tf\setminus\{\epsilon_i^2-\epsilon_i\mid i\in Sp\})\cup R^{\Tf_1}=I\cup R^{\Tf_1},$$ where $R^{\Tf_1}=\{\rho_{i}\rho_{i^*}\mid i\in Sp\}$({\em cf.} Figure~\ref{fig: rho_irho_i*=0}), which implies $A_{\Tf^*}\cong KQ^*/I^*$.
  \end{proof}

\begin{remark}\label{r: bijection between e and arcs} 
       By \cite[Remarks 1.4 and 3.8]{HZZ}, there is a bijection  between the set
\[
\{\varepsilon_j\mid j\in Q_0{\setminus}Sp\}\cup\{\varepsilon_i-\epsilon_i,\,\epsilon_i\mid i\in Sp\}
\]
and $\Tf^{\bowtie}$ given by:
\begin{itemize}
 \item mapping $\varepsilon_j$ to the arc $\gamma$ in $\Tf\setminus\Tf^{sp}$  such that $j$ is indexed by $\gamma$;
  \item mapping $\varepsilon_{i}-\epsilon_i$ (resp. $\epsilon_i$) to the tagged arc $\ba_P^-$ (resp. $\ba_P^+$) at the puncture $P$ such that $i$ is indexed by $\lf_P$.
\end{itemize}
Moreover, because  $A^{sg}$ and $A^{\bowtie}$  are Morita equivalent and the set
\[
S^{\bowtie}=\{\varepsilon_j\mid j\in Q_0{\setminus}Sp\}\cup\{\varepsilon_{i^-},\varepsilon_{i^+}\mid i\in Sp\}
\]
forms a complete set of primitive orthogonal idempotents of $A^{{\bowtie}}$, so
the Morita equivalence then induces a bijection from  $S^{\bowtie}$ to $\Tf^{\bowtie}$
given by
\begin{itemize}
  \item mapping $\varepsilon_j$ to the arc $\gamma$ in $\Tf\setminus\Tf^{sp}$  such that $j$ is indexed by $\gamma$;
  \item mapping $\varepsilon_{i^-}$ (resp. $\varepsilon_{i^+}$) to the tagged arc $\ba_P^-$ (resp. $\ba_P^+$) at the puncture $P$ such that $i$ is indexed by $\lf_P$.
\end{itemize}

\end{remark}

 \section{Intersection vectors between skew-tiling and unfolded tiling}\label{s:s3}
\setcounter{equation}{0}
Let $(\Sf,\Mf,\Pf, \Tf)$ be a  skew-tiling,  $\Tf^{\bowtie}$ the tagged version of $\Tf$, and $(\Sf^*,\Mf^*,\Tf^*)$ the associated {unfolded tiling}.
\subsection{Tagged permissible arcs} 
We use the following notion of permissible arcs from \cite[Definition~2.2]{HZZ}. 
\begin{definition}\label{def:permissible}
	An arc $\gamma$ on $\Sf$ is called \emph{permissible} with respect to $\Tf$ if it satisfies the following conditions:
	\begin{enumerate}
		\item[(P1)] The starting/ending segment of $\gamma$ matches one of the local configurations shown in Figure~\ref{fig:ro}, where the blue arcs belong to $\Tf$ or are boundary segments.
		\item[(P2)] Whenever $\gamma$ consecutively crosses two non-boundary edges $x$ and $y$ (which may coincide) of a tile $\triangle$, then $x$ and $y$ are adjacent sides incident to $p_{\eta}\in\Mf$ and $\gamma$ cuts off an angle from $\triangle$, as depicted in Figure~\ref{fig:local tri}.
	\end{enumerate}
    \end{definition}
	We denote by  $\PAS$ the set of permissible arcs with respect to $\Tf$.

A {\it tagged permissible} arc on $\Sf$ is a tagged arc whose untagged version is permissible  with respect to $\Tf$. We denote the set of all such tagged permissible arcs by $\PTAS$. Furthermore, let $\mathscr{R}(\Sf)$ be the set of all finite multisets consisting of pairwise compatible tagged permissible arcs. By (P1), we know for any arc $\af\in \Tf^{\bowtie}$, $\af\notin \PTAS$. Moreover, we know for a a pair $\{\gamma,\gamma'\}$ of conjugate arcs, if one of $\{\gamma,\gamma'\}$ is in $\Tf$, the other one must lie in $\Tf$. Thus for any arc $\gamma\in \PTAS$, 
\begin{equation*} 
\begin{aligned}
    \Int_{\Tf^{\bowtie}}(\af|\gamma)=\Int_{\Tf^{\bowtie}}^{A}(\af|\gamma)+\Int_{\Tf^{\bowtie}}^D(\af|\gamma).
    \end{aligned}
\end{equation*}
    
	\begin{figure}[htpb]
		\begin{tikzpicture}[scale=1.25]
			\clip(-1.2,-1.2) rectangle (1.2,1.2);
			\draw[blue, thick, bend left=10](0,1)to(0,-1);
			\draw[blue, thick](0,1)to(1,0);
			\draw[blue,thick,dashed,bend right=10](0,-1)to(1,0);
			\draw[red, thick, bend right=10](-1,0)to(1,0);
			\draw[red,thick](-1,0)node[below]{$\gamma$};
			\draw[thick](0,1)node{$\bullet$}(0,-1)node{$\bullet$}(1,0)node{$\bullet$};
		\end{tikzpicture}\qquad
		\begin{tikzpicture}[scale=1.25]
			\clip(-1.2,-1.2) rectangle (1.2,1.2);
			\draw[blue,thick,bend right=60](0,1)to(0,-1);
			\draw[blue,thick,bend left=60](0,1)to(0,-1);
			\draw[ultra thick, fill=gray!20](0,0)circle(.1);
			\draw[red,thick,smooth](-1,.4)to[out=-5,in=170](0,.3)to[out=-10,in=90](.3,0)to[out=-90,in=70](0,-1);
			\draw[red](-1,0.1)node{$\gamma$};
			\draw[thick](0,1)node{$\bullet$}(0,-1)node{$\bullet$};
 		\end{tikzpicture}\qquad
		\begin{tikzpicture}[scale=1.25]
			\clip(-1.2,-1.2) rectangle (1.2,1.2);
			\draw[blue,thick](0,1)to[out=-135,in=80](-.6,0)to[out=-100,in=180](0,-1)to[out=0,in=-80](.6,0)to[out=100,in=-45](0,1);
			\draw[red,thick](0,0)to(0,-1.2);
			\draw[red] (-.1,-.4)node[right]{$\gamma$};
			\draw(0,1)node{$\bullet$};
			\draw(0,0)node{$\bullet$};
		\end{tikzpicture}
		\caption{Condition (P1)}
		\label{fig:ro}
	\end{figure}
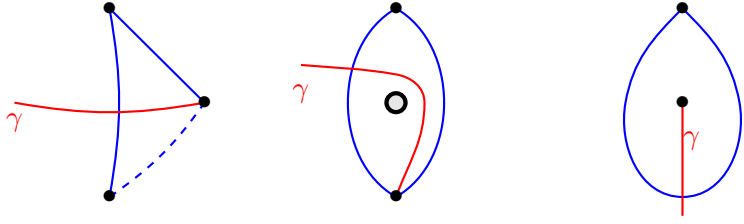
	\begin{figure}[htpb]\centering
		\begin{tikzpicture}[xscale=0.8,yscale=0.8]
				\draw[thick,fill=black!20] (9,0.17)arc (250:290:3);
				\node at (10,0.2){$p_{\eta}$};
				\fill(10,0) circle(1.5pt);
				\draw[thick,color=blue] (10,0) .. controls (8.6,-0.9) and (8.4,-1.1) .. (8,-2.2);
				\draw[thick,color=blue,dashed] (8,-2.2) -- (7.8,-3);
				\draw[thick,color=blue] (10,0) .. controls (10.4,-0.9) and (10.6,-1.1) .. (11,-2.2);
				\draw[thick,color=blue,dashed] (11,-2.2) arc (20:-40:1);
				\draw[thick,color=red] (8,-1) .. controls (9,-1.5) and (10,-1) .. (12,-1);
				\draw[thick,color=red,dashed] (12,-1) arc (90:50:1.5);
				\draw[thick,color=red,dashed] (8,-1) -- (7,-0.5);
				\node at (8.5,-.5){{\color{blue}$x$}};
				\node at (10.5, -.5){{\color{blue}$y$}};
                \node at (11, -1.5){{\color{red}$\gamma$}};
				\draw (9.7,-.8)node{$\triangle$};
			\end{tikzpicture}
		\caption{Condition (P2)}
		\label{fig:local tri}
	\end{figure}
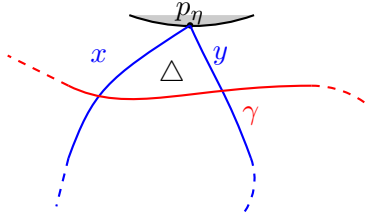


\subsection{Intersection vectors of tagged permissible arcs}
 For a tagged permissible arc $\af$ and a multiset $\Mcal$ of tagged permissible arcs, we define their \emph{intersection number} as: 
\[\Int(\af|\Mcal)=\Int(\Mcal|\af):=\sum_{\gamma\in \Mcal}\Int(\af| \gamma).\]
Given a tagged permissible arc $\gamma$, its {\it intersection vector} with respect to $\Tf^{\bowtie}$ is defined as:
\[\underline{\Int}_{\Tf^{\bowtie}}(\gamma)=(\Int(\af|\gamma))_{\af\in\Tf^{\bowtie} }.\]
For a multiset $\Mcal$ in $\mathscr{R}(\Sf)$, the {\it intersection vector} of $\mathcal{M}$ with respect to $\Tf^{\bowtie}$ is given by:  
\[\underline{\Int}_{\Tf^{\bowtie}}(\mathcal{M}):=\sum_{\gamma\in\mathcal{M}}\underline{\Int}_{\Tf^{\bowtie}}(\gamma).\]
Note that $\underline{\Int}_{\Tf^{\bowtie}}(\mathcal{M})=(\Int(\mathbf{a}|\mathcal{M}))_{\mathbf{a}\in\Tf^{\bowtie}}.
$\vspace{0.15cm}
\subsection{Permissible arcs from skew-tiling to unfolded tiling}
In this subsection, we construct a map $\phi$ from $\PTAS$ to $\mathbb{PA}(\mathbf{S}^*)$ and a map $\Phi$ from $\mathscr{R}(\Sf)$ to $\mathscr{R}(\Sf^*)$.

\begin{figure}[ht]
\begin{tikzpicture}[xscale=0.8,yscale=0.8]
\draw[thick,color=blue] (0,-2)\nn .. controls (-0.5,-1.5) and (-0.5,-0.5) .. (0,0)\nn;	
\draw[thick,color=blue] (0,-2) .. controls (0.5,-1.5) and (0.5,-0.5) .. (0,0);	
\draw[thick,blue](-0.3,-1)node[left]{$\af_P^{-}$};
\draw[thick,blue](0,0)node[above]{P}(0.3,-1)node[right]{$\af_P^{+}$}(0,-2)node[below]{$b_P$}(1,0)(0,0.8);

	\draw[red, thick, bend right=10](0,0)to(2,0);
			\draw[red,thick](2,0)node[below]{$\gamma$};
		\end{tikzpicture}
  \begin{tikzpicture}[xscale=0.8,yscale=0.8]
\draw[white,thick,=>] (0,2) to (2,2);
\draw[red,thick,->] (0,0) to (2,0);
\draw[white,thick,->] (0,-2) to (2,-2);
	\end{tikzpicture}
\begin{tikzpicture}[xscale=0.8,yscale=0.8]
\draw[thick,fill=black!20] (0,0.45) circle(12pt);
\draw (0,-2)\nn  ;	
\draw[blue,thick,>=stealth] (0,-2)\nn to (0,0)\nn;
\draw[thick,color=blue] (0,-2) .. controls (-2,-1.5) and (-2,2) .. (0,2);	
\draw[thick,color=blue] (0,-2) .. controls (2,-1.5) and (2,2) .. (0,2);	
\draw[thick,blue](0,-0.5)node[right]{$\af_{P}^{-}$} (0,2)node[above]{$\lf_P$}(0,0)(0,-2)node[below] {$b_P$}(0,0.8)node[above]{$B_P$};
	\draw[red, thick, bend right=10](0,0)to(2,0);
    \draw[red,thick](2,0)node[below]{$\phi(\gamma)$};
\end{tikzpicture}
\caption{From $\gamma$ with some puncture as one endpoint  to  $\phi(\gamma)$}
\label{f:tagged arcs to loop}
\end{figure}
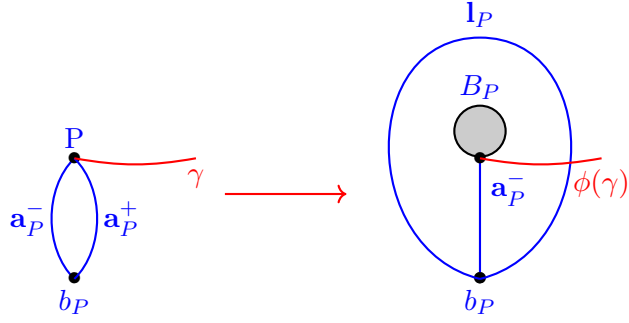

For a tagged permissible arc $\gamma\in\PTAS$, let $\gamma^{\circ}$ denote its untagged version. We define the map $\phi$ as follows:
\begin{itemize}
    \item If neither endpoint of $\gamma$ is a puncture, then after transforming $(\Sf,\Mf,\Pf,\Tf)$ into $(\Sf^*,\Mf^*,\Tf^*)$, we define $\phi(\gamma)$ to be the arc $\gamma$ itself, viewed as an arc in $(\Sf^*,\Mf^*,\Tf^*)$.
    \item If one endpoint of $\gamma$ is a puncture $P$, then after transforming $(\Sf,\Mf,\Pf,\Tf)$ into $(\Sf^*,\Mf^*,\Tf^*)$ and placing $B_P$ as  shown in Figure~\ref{f:tagged arcs to loop}, we define $\phi(\gamma)$ to be the arc  $\gamma^{\circ}$, again viewed as an arc in the new marked surface. 
\end{itemize}
In both cases, $\phi(\gamma)$ is given by the untagged arc $\gamma^{\circ}$ considered within the new marked surface $(\Sf^*,\Mf^*,\Tf^*)$. Furthermore, $\phi(\gamma)=\gamma^\circ$ is also permissible with respect to $T^*$, hence $\phi(\gamma)\in \mathbb{PA}(\mathbf{S}^*)$. 


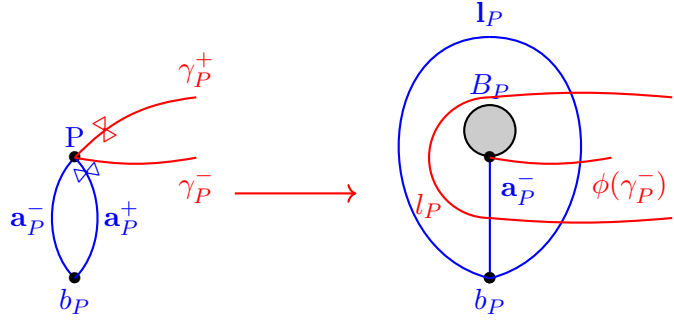
\begin{figure}[ht]
\begin{tikzpicture}[xscale=0.8,yscale=0.8]
\draw[thick,color=blue] (0,-2)\nn .. controls (-0.5,-1.5) and (-0.5,-0.5) .. (0,0)\nn;	
\draw[thick,color=blue] (0,-2) .. controls (0.5,-1.5) and (0.5,-0.5) .. (0,0);	
\draw[blue](0.2,-.25)node[rotate=20]{$\bowtie$};	
\draw[thick,blue](-0.3,-1)node[left]{$\af_P^{-}$};
\draw[thick,blue](0,0)node[above]{P}(0.3,-1)node[right]{$\af_P^{+}$}(0,-2)node[below]{$b_P$}(1,0)(0,0.8);

\draw[red, thick, bend left=20](0,0)to(2,1);
\draw[red](0.5,0.47)node[rotate=120]{$\bowtie$};	
		\draw[red,thick](2,1)node[above]{$\gamma_P^{+}$};
	\draw[red, thick, bend right=10](0,0)to(2,0);
			\draw[red,thick](2,0)node[below]{$\gamma_{P}^{-}$};
		\end{tikzpicture}
  \begin{tikzpicture}[xscale=0.8,yscale=0.8]
\draw[white,thick,=>] (0,2) to (2,2);
\draw[red,thick,->] (0,0) to (2,0);
\draw[white,thick,->] (0,-2) to (2,-2);
	\end{tikzpicture}
\begin{tikzpicture}[xscale=0.8,yscale=0.8]
\draw[thick,fill=black!20] (0,0.45) circle(12pt);
\draw (0,-2)\nn  ;	
\draw[blue,thick,>=stealth] (0,-2)\nn to (0,0)\nn;
\draw[thick,color=blue] (0,-2) .. controls (-2,-1.5) and (-2,2) .. (0,2);	
\draw[thick,color=blue] (0,-2) .. controls (2,-1.5) and (2,2) .. (0,2);	
\draw[thick,blue](0,-0.5)node[right]{$\af_{P}^{-}$} (0,2)node[above]{$\lf_P$}(0,0)(0,-2)node[below] {$b_P$}(1,0)(0,0.8)(0,0.8)node[above]{$B_P$};
\draw[red,thick,smooth](3,1)to[bend right=5](0,1)to
[out=-180,in=90](-1,0)to [out=-90,in=180] (0,-1)to [bend right=5](3,-1);
			\draw[red](-1,-0.8)node{$l_{P}$};

\draw[red, thick, bend right=10](0,0)to(2,0);
			\draw[red,thick](2.3,0)node[below]{$\phi(\gamma_{P}^{-})$};

\end{tikzpicture}
\caption{From the conjugate pair $\{\gamma_{P}^{-},\gamma_P^{+}\}$ to the loop $l_{P}$}
\label{f:pairs to loop}
\end{figure}

 Let $\Mcal\in\mathscr{R}(\Sf)$. For any puncture $P$, let $\{\gamma_{P}^{-}, \gamma_P^{+}\}$ denote the pair of conjugate arcs in $\Mcal$ tagged at $P$,  where $\gamma_{P}^{-}$ is the arc tagged plain at the end incident to $P$ and $\gamma_P^{+}$ is the arc tagged notched at the end incident to $P$. Let $l_{P}$ be the unique corresponding permissible loop in $(\Sf^*,\Mf^*,\Tf^*)$ that closely wraps around $\phi(\gamma_{P}^{-})$ and $B_P$, as shown in Figure \ref{f:pairs to loop}.
Denote by $\Mcal_P$ be the sub-multiset of $\Mcal$ consisting of  $\min(\Int_{\Tf^{\bowtie}}^D(\af_P^{-}|\Mcal),\Int_{\Tf^{\bowtie}}^D(\af_P^{+}|\Mcal))$ copies of conjugate pairs $\{\gamma_{P}^{-}, \gamma_P^{+}\}$ in  $\Mcal$. Denote by 
$$\Mcal_1=\bigcup_{P\in\Pf}\Mcal_P, \ \ \ \ \Mcal_2=\Mcal\setminus \Mcal_1.$$
Then there is no pair of conjugate arcs in $\Mcal_2$.
We define $\Mcal^*$ as the multiset of  arcs in $(\Sf^*,\Mf^*,\Tf^*)$ obtained from $\Mcal$ via the following procedure:
\begin{itemize}
    \item [$\textbf{Step } 1$:] For each puncture $P$, replace each pair of conjugate arcs $\{\gamma_{P}^{-}, \gamma_P^{+}\}$ in  $\Mcal_P$  with the unique corresponding permissible loop $l_{P}$;
    \item [$\textbf{Step } 2$:] Replace each arc $\gamma$ in $\Mcal_2$  with its image $\phi(\gamma)$, \ie $\gamma^{\circ}$.
\end{itemize}
It is clear that  $\Mcal^*\in \mathscr{R}(\Sf^*)$. This defines a map $\Phi: \mathscr{R}(\Sf)\rightarrow \mathscr{R}(\Sf^*)$
by setting $\Phi(\Mcal)=\Mcal^*$ for each $\Mcal\in \mathscr{R}(\Sf)$.

\subsection{Permissible arcs from  unfolded tiling to skew-tiling}\label{ss:Permissible arcs from  unfolded tiling to skew-tiling}
 Let $\Mcal\in\mathscr{R}(\Sf)$.
Denote by
  $$S_{\Mcal}=\left\{P\in \Pf\Big| \Int_{\Tf^{\bowtie}}^D(\af_P^{+}|\Mcal)-\Int_{\Tf^{\bowtie}}^D(\af_P^{-}|\Mcal)< 0\right\}.$$ 
Then  $\Mcal$ can be recovered from $\Mcal^{*}$ and  $S_{\Mcal}$ through the following process:
\begin{itemize}
\item [$\textbf{Step}\  1^{\circ}$:] For each puncture $P$ with $P\in S_{\Mcal}$,  tag  notched at the end of each arc in  $\Mcal^{*}$ which is  incident to  $P$. For any arc $\gamma\in\mathcal M^{*}$, write $\hat\gamma$ for the arc produced from $\gamma$ and  write $\hat{\Mcal}$ for the multiset produced from $\Mcal^{*}$ after this step.
\item  [$\textbf{Step}\  2^{\circ}$:] For every puncture $P\in\Pf$, replace each  once-punctured loop $\hat{l}_P$ in $\hat\Mcal^{*}$ enclosing $P$ with the corresponding conjugate pair
$\{\gamma_{P}^{-},\,\gamma_{P}^{+}\}$, where   $\gamma_P^{-}$ and $\gamma_P^{+}$ are tagged in the same way with $\hat{l}_P$ at the ends incident to the base point of $\hat{l}_P$.
\end{itemize}

\subsection{Intersection vectors from skew-tiling to unfolded tiling}
Recall that $$\Tf^*=(\Tf^{\bowtie}\setminus \{\af_P^+|P\in\Pf\})\bigcup \{\lf_P\in \Tf|P\in\Pf\}.$$ Then we have

\begin{lemma}\label{l:skew-gentle to gentle}
Let $\Mcal\in\mathscr{R}(\Sf)$. For any puncture $P\in\Pf$, we have
 \begin{itemize}
 \item[(1)]  If $P\notin S_{\Mcal}$, then $\Int_{\Tf^*}(\af_P^{-}|\Mcal^*)= \Int_{\Tf^{\bowtie}}(\af_{P}^{-}|\Mcal).$\vspace{0.15cm}
 \item[(2)] If $P\in S_{\Mcal}$, then $\Int_{\Tf^*}(\af_P^{-}|\Mcal^*)=\Int_{\Tf^{\bowtie}}(\af_{P}^{+}|\Mcal).$\vspace{0.15cm}
\item[(3)] In either case, $\Int_{\Tf^*}(\lf_P|\Mcal^*)=\Int_{\Tf^{\bowtie}}(\af_{P}^{-}|\Mcal)+\Int_{\Tf^{\bowtie}}(\af_{P}^{+}|\Mcal).$
  \end{itemize}
 On the other hand, for $\af\in \Tf^{\bowtie}$ such that  no puncture is the endpoint of $\af$, then 
  \begin{itemize}
      \item [(4)] $\Int_{\Tf^*}(\af|\Mcal^*)= \Int_{\Tf^{\bowtie}}(\af|\Mcal).$
  \end{itemize}
\end{lemma}
\begin{proof} Note that for any arc $\gamma\in \Mcal$, we have
$\Int_{\Tf^{\bowtie}}^{A}(\af_{P}^{-}|\gamma)=\Int_{\Tf^{\bowtie}}^{A}(\af_{P}^{+}|\gamma)$, and 
\begin{equation*} 
\begin{aligned}
    \Int_{\Tf^{\bowtie}}(\af|\gamma)=\Int_{\Tf^{\bowtie}}^{A}(\af|\gamma)+\Int_{\Tf^{\bowtie}}^D(\af|\gamma)
    \end{aligned}
\end{equation*}
for any $\af\in \Tf^{\bowtie}$.
Since the unfolded tiling $(\Sf^*,\Mf^*,\Tf^*)$ contains no puncture,  we have $$\Int_{\Tf^*}(\af|\beta)=\Int_{\Tf^*}^{A}(\alpha|\beta)$$ for any  $\alpha\in \Tf^*$ and  $\beta \in\mathbb{PA}(\Sf^*)$.
Recall that for each puncture $Q$, $\Mcal_Q$ is the sub-multiset of $\Mcal$ consisting of  $\min(\Int_{\Tf^{\bowtie}}^D(\af_Q^{-}|\Mcal),\Int_{\Tf^{\bowtie}}^D(\af_Q^{+}|\Mcal))$ copies of conjugate pairs $\{\gamma_{Q}^{-}, \gamma_Q^{+}\}$ in  $\Mcal$ and 
$$\Mcal_1=\bigcup_{Q\in\Pf}\Mcal_Q,\ \ \  \Mcal_2=\Mcal\setminus \Mcal_1.$$
It is clear that  there is no pair of conjugate arcs in $\Mcal_2$,
each arc in $\Mcal_1^*$ is a loop $l_Q$ derived from some conjugate pair $\{\gamma_Q^-,\gamma_Q^+\}$ in $\Mcal_1$ and
each arc in $\Mcal_2^*$ is the untagged version $\gamma^{\circ}$ of some arc $\gamma\in\Mcal_2$. 

For  any   $\gamma\in\Mcal_2$, we have
\begin{align}\label{g:a-|gamma0}
     \Int_{\Tf^*}(\af_{P}^{-}|\gamma^{\circ})=\Int_{\Tf^*}^{A}(\af_{P}^{-}|\gamma^{\circ})=\Int_{\Tf^{\bowtie}}^{A}(\af_{P}^{-}|\gamma).
\end{align}
For any conjugate pair $\{\gamma_{Q}^{-},\gamma_{Q}^{+}\}$ in $\Mcal_1$, since $(\gamma_{Q}^{-})^{\circ}$ can be closely wrapped by $l_{Q}$,  we can regard $l_{Q}$ as  $((\gamma_{Q}^{-})^{\circ})^{(2)}$ when $l_Q$ passes through $ \af_{P}^{-}$. Thus if $Q\neq P$, then
\begin{align}\label{g: a-|l_Q}
   \Int_{\Tf^*}(\af_{P}^{-}|l_{Q})
   =2\Int_{\Tf^*}^{A}(\af_{P}^{-}|(\gamma^-_{Q})^\circ)
   =&2\Int_{\Tf^{\bowtie}}^{A}(\af_{P}^{-}|\gamma^-_{Q})\\
     =&\Int_{\Tf^{\bowtie}}^{A}(\af_{P}^{-}|\gamma^-_{Q})+\Int_{\Tf^{\bowtie}}^{A}(\af_{P}^{-}|\gamma^+_{Q}) \notag
\end{align}
({\em cf.} $\Int(\af_{P_2}^-|l_{P_3})$ in Figure \ref{fig:intersection of  l_{p_3} and T}).  Else if $Q=P$, then
\begin{align}\label{g:a-|l_p}
      \Int_{\Tf^*}(\af_{P}^{-}|l_{P})
   =1+2\Int_{\Tf^*}^{A}(\af_{P}^{-}|(\gamma^-_{P})^\circ)
   =&1+2\Int_{\Tf^{\bowtie}}^{A}(\af_{P}^{-}|\gamma^-_{P})\\
    =&1+\Int_{\Tf^{\bowtie}}^{A}(\af_{P}^{-}|\gamma^-_{P})+\Int_{\Tf^{\bowtie}}^{A}(\af_{P}^{-}|\gamma^+_{P})\notag
\end{align}
({\em cf.} $\Int(\af_{P_3}^-|l_{P_3})$ in Figure \ref{fig:intersection of  l_{p_3} and T}).

 For $(1)$, when $P\notin S_{\Mcal}$,  we have $\Int_{\Tf^{\bowtie}}^D(\af_P^{+}|\Mcal)-\Int_{\Tf^{\bowtie}}^D(\af_P^{-}|\Mcal)\ge 0$. If $\Int_{\Tf^{\bowtie}}^D(\af_P^{-}|\Mcal)=0$, then $\Mcal_P=\emptyset$ and all tagged arcs in $\Mcal$ with $P$ as an endpoint  are tagged plain at the end incident to $P$. Else if $\Int_{\Tf^{\bowtie}}^D(\af_P^{-}|\Mcal)>0$,  since $\Int_{\Tf^{\bowtie}}^D(\af_P^{+}|\Mcal)\ge \Int_{\Tf^{\bowtie}}^D(\af_P^{-}|\Mcal)>0$ and 
any two arcs in $\Mcal$ are mutually compatible, thus any  tagged arc $\gamma\in \Mcal\setminus \Mcal_P$ with $P$ as an endpoint must be $\gamma_{P}^{-}$. 
Hence, in either case, we  have  $$\Int_{\Tf^{\bowtie}}^{D}(\af_{P}^{-}|\gamma)=0 \text{ for any }\gamma\in \Mcal\setminus\Mcal_P.$$ 
Then  for  any   $\gamma\in\Mcal_2$, by~(\ref{g:a-|gamma0}), we have
\begin{align*}
     \Int_{\Tf^*}(\af_{P}^{-}|\gamma^{\circ})=\Int_{\Tf^{\bowtie}}^{A}(\af_{P}^{-}|\gamma)=\Int^A_{\Tf^{\bowtie}}(\af_{P}^{-}|\gamma)+\Int_{\Tf^{\bowtie}}^{D}(\af_{P}^{-}|\gamma)=\Int_{\Tf^{\bowtie}}(\af_{P}^{-}|\gamma).
\end{align*}
For any conjugate pair $\{\gamma_{Q}^{-},\gamma_{Q}^{+}\}$ in $\Mcal_1$, if $Q\neq P$, by~(\ref{g: a-|l_Q})
\begin{align*}
   \Int_{\Tf^*}(\af_{P}^{-}|l_{Q})
     =&\Int_{\Tf^{\bowtie}}^{A}(\af_{P}^{-}|\gamma^-_{Q})+\Int_{\Tf^{\bowtie}}^{A}(\af_{P}^{-}|\gamma^+_{Q}) \\
     =&\Int^A_{\Tf^{\bowtie}}(\af_{P}^{-}|\gamma^-_{Q})+\Int^A_{\Tf^{\bowtie}}(\af_{P}^{-}|\gamma^+_{Q})+(\Int_{\Tf^{\bowtie}}^{D}(\af_{P}^{-}|\gamma^-_{Q})+\Int_{\Tf^{\bowtie}}^{D}(\af_{P}^{-}|\gamma^+_{Q}))\\
     =&\Int_{\Tf^{\bowtie}}(\af_{P}^{-}|\gamma^-_{Q})+\Int_{\Tf^{\bowtie}}(\af_{P}^{-}|\gamma^+_{Q}).
\end{align*}
  Else if $Q=P$, by~(\ref{g:a-|l_p}), we have
\begin{align*}
      \Int_{\Tf^*}(\af_{P}^{-}|l_{P})
    =&1+\Int_{\Tf^{\bowtie}}^{A}(\af_{P}^{-}|\gamma^-_{P})+\Int_{\Tf^{\bowtie}}^{A}(\af_{P}^{-}|\gamma^+_{P})\\
    =&(\Int_{\Tf^{\bowtie}}^{D}(\af_{P}^{-}|\gamma^-_{P})+\Int_{\Tf^{\bowtie}}^{D}(\af_{P}^{-}|\gamma^+_{P}))+\Int_{\Tf^{\bowtie}}^{A}(\af_{P}^{-}|\gamma^-_{P})+\Int_{\Tf^{\bowtie}}^{A}(\af_{P}^{-}|\gamma^+_{P})\\
    =&\Int_{\Tf^{\bowtie}}^{}(\af_{P}^{-}|\gamma^-_{P})+\Int_{\Tf^{\bowtie}}^{}(\af_{P}^{-}|\gamma^+_{P})
\end{align*}
({\em cf.} $\Int(\af_{P_3}^-|l_{P_3})$ in Figure \ref{fig:intersection of  l_{p_3} and T}). In total, statement $(1)$ holds.
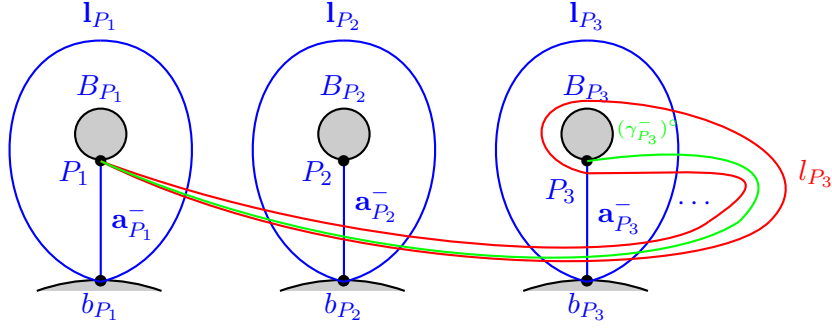
\begin{figure}
    \centering
    \begin{tikzpicture}[xscale=0.8,yscale=0.8]
 \draw[thick,fill=black!20] (1,-2.15)arc (70:110:3);
\draw[thick,fill=black!20] (0,0.45) circle(12pt);
\draw (0,-2)\nn  ;	
\draw[blue,thick,>=stealth] (0,-2)\nn to (0,0)\nn;
\draw[thick,color=blue] (0,-2) .. controls (-2,-1.5) and (-2,2) .. (0,2);	
\draw[thick,color=blue] (0,-2) .. controls (2,-1.5) and (2,2) .. (0,2);	
\draw[thick,blue](0,0)(0,-2)node[below] {$b_{P_1}$}(0,0.8)node[above]{$B_{P_1}$};
\draw[thick,fill=black!20](5,-2.15) arc(70:110:3);
\draw[thick,fill=black!20] (4,0.45) circle(12pt);
\draw (4,-2)\nn;	
\draw[blue,thick,>=stealth] (4,-2)\nn to (4,0)\nn;
\draw[thick,color=blue] (4,-2) .. controls (2,-1.5) and (2,2) .. (4,2);	
\draw[thick,color=blue] (4,-2) .. controls (6,-1.5) and (6,2) .. (4,2);	
\draw[thick,blue](4,0)(4,-2)node[below] {$b_{P_2}$}(4,0.8)node[above]{$B_{P_2}$};
\draw[thick,fill=black!20](9,-2.15) arc(70:110:3);
\draw[thick,fill=black!20] (8,0.45) circle(12pt);
\draw (8,-2)\nn;	
\draw[blue,thick,>=stealth] (8,-2)\nn to (8,0)\nn;
\draw[thick,color=blue] (8,-2) .. controls (6,-1.5) and (6,2) .. (8,2);	
\draw[thick,color=blue] (8,-2) .. controls (10,-1.5) and (10,2) .. (8,2);	
\draw[thick,blue](8,0)(8,-2)node[below] {$b_{P_3}$}(8,0.8)node[above]{$B_{P_3}$};
\draw[thick,blue](0,2)node[above]{$\lf_{P_1}$};
\draw[thick,blue](4,2)node[above]{$\lf_{P_2}$};
\draw[thick,blue](8,2)node[above]{$\lf_{P_3}$};

\draw[blue,thick] (9.8,-0.7) node{$\cdots$};
\draw[thick,blue](0,-.2)node[left] {${P_1}$}(4,-.2)node[left] {${P_2}$}(8,-.5)node[left] {${P_3}$};

\draw[thick,blue](0,-1)node[right] {$\af^{-}_{P_1}$}(4,-.7)node[right] {$\af^{-}_{P_2}$}(8,-.9)node[right] {$\af^{-}_{P_3}$};

\draw[thick,color=red] (0,0) .. controls (4,-2) and (10,-2) .. (11,-1);
\draw[thick,color=red] (8,1) .. controls (10,1) and (12,0) .. (11,-1);

\draw[thick,color=red] (0,0) .. controls (4,-1.5) and (9,-1.8) .. (10,-1);
\draw[thick,color=red] (8,-0.2) .. controls (10,-0.2) and (11.5,0) .. (10,-1);

\draw[thick,color=red] (8,-0.2) .. controls (7,0) and (7,1) .. (8,1);

\draw[thick,color=green] (8,0) .. controls (10,0.3) and (11.5,0) .. (10.5,-1);
\draw[thick,color=green] (0,0) .. controls (4,-1.8) and (9,-2) .. (10.5,-1);
\draw[thick,red](11.3,-.2)node[right] {$l_{P_3}$};
\draw[thick,green](9,0.1)node[above] {\tiny $(\gamma^-_{P_3})^{\circ}$};

\end{tikzpicture}
    \caption{Intersections between $l_{P_3}$ and $\Tf^*$}
    \label{fig:intersection of  l_{p_3} and T}
\end{figure}

\vspace{0.3cm}

For $(2)$, by applying a similar argument as in $(1)$, we obtain the desired result.

\vspace{0.3cm}

For $(3)$, since $\lf_P$  closely wraps  $\af_P^{-}$ and $B_P$ in the tiling $(\Sf^*,\Mf^*,\Tf^*)$, it can be viewed as $(\af_P^{-})^{(2)}$ whenever $\af_P^-$ passes through  another arc.
Then for each   $\gamma\in\Mcal_2$, we consider the following cases:
\begin{itemize}
    \item [Case 1.1:] If neither endpoint of $\gamma$ is a puncture, then 
    $\Int_{\Tf^{\bowtie}}^{D}(\af_{P}^{-}|\gamma)+\Int_{\Tf^{\bowtie}}^{D}(\af_{P}^{+}|\gamma)=0$, and
 $$\Int_{\Tf^*}(\lf_P|\gamma^{\circ})=\Int^A_{\Tf^*}(\lf_P|\gamma^{\circ})=2\Int^A_{\Tf^*}(\af_{P}^{-}|\gamma^{\circ})=\Int_{\Tf^{\bowtie}}^{A}(\af_{P}^{-}|\gamma)+\Int_{\Tf^{\bowtie}}^{A}(\af_{P}^{+}|\gamma).$$ 
Thus, $\Int_{\Tf^*}(\lf_P|\gamma^{\circ})=\Int_{\Tf^{\bowtie}}(\af_{P}^{-}|\gamma)+\Int_{\Tf^{\bowtie}}(\af_{P}^{+}|\gamma).$

\vspace{0.2cm}
    \item [Case 1.2:] If one endpoint of $\gamma$ is the puncture $P$ and the other endpoint is not, then $\Int_{\Tf^{\bowtie}}^{D}(\af_{P}^{-}|\gamma)+\Int_{\Tf^{\bowtie}}^{D}(\af_{P}^{+}|\gamma)=1$, and 
     $$\Int_{\Tf^*}(\lf_P|\gamma^{\circ})=\Int^A_{\Tf^*}(\lf_P|\gamma^{\circ})=1+2\Int^A_{\Tf^*}(\af_{P}^{-}|\gamma^{\circ})=1+\Int_{\Tf^{\bowtie}}^{A}(\af_{P}^{-}|\gamma)+\Int_{\Tf^{\bowtie}}^{A}(\af_{P}^{+}|\gamma)$$  
   ({\em cf.} $\Int(\lf_{P_3}|(\gamma^-_{P_3})^{\circ})$ in Figure \ref{fig:intersection of  l_{p_3} and T}).  Hence, $\Int_{\Tf^*}(\lf_P|\gamma^{\circ})=\Int_{\Tf^{\bowtie}}(\af_{P}^{-}|\gamma)+\Int_{\Tf^{\bowtie}}(\af_{P}^{+}|\gamma).$

\vspace{0.2cm}
    \item [Case 1.3:]  If both endpoints of $\gamma$ are the puncture $P$, then $\Int_{\Tf^{\bowtie}}^{D}(\af_{P}^{-}|\gamma)+\Int_{\Tf^{\bowtie}}^{D}(\af_{P}^{+}|\gamma)=2$, and 
     $$\Int_{\Tf^*}(\lf_P|\gamma^{\circ})=\Int^A_{\Tf^*}(\lf_P|\gamma^{\circ})=2+2\Int^A_{\Tf^*}(\af_{P}^{-}|\gamma^{\circ})=2+\Int_{\Tf^{\bowtie}}^{A}(\af_{P}^{-}|\gamma)+\Int_{\Tf^{\bowtie}}^{A}(\af_{P}^{+}|\gamma).$$  
     Thus, $\Int_{\Tf^*}(\lf_P|\gamma^{\circ})=\Int_{\Tf^{\bowtie}}(\af_{P}^{-}|\gamma)+\Int_{\Tf^{\bowtie}}(\af_{P}^{+}|\gamma).$

\end{itemize}

Therefore, for any $\gamma\in\Mcal_2$,  we have 
\begin{eqnarray}\label{g: case wilde{a_p}+ with gamma_0}
    \Int_{\Tf^*}(\lf_P|\gamma^{\circ})=\Int_{\Tf^{\bowtie}}(\af_{P}^{-}|\gamma)+\Int_{\Tf^{\bowtie}}(\af_{P}^{+}|\gamma).
\end{eqnarray}

For each  conjugate pair $\{\gamma_{Q}^{-},\gamma_{Q}^{+}\}$ in $\Mcal_1$, if $Q\neq P$, we have the following cases:
\begin{itemize}
    \item [Case 2.1]  $P$ is not the  basepoint of $l_Q$, then $P$ is not the endpoint of either $\gamma_Q^-$ or $\gamma_Q^+$.
   Hence, $$\Int_{\Tf^{\bowtie}}^{D}(\af_{P}^{-}|\gamma_{Q}^{-})+\Int_{\Tf^{\bowtie}}^{D}(\af_{P}^{-}|\gamma_{Q}^{+})+\Int_{\Tf^{\bowtie}}^{D}(\af_{P}^{+}|\gamma_{Q}^{-})+\Int_{\Tf^{\bowtie}}^{D}(\af_{P}^{+}|\gamma_{Q}^{+})=0.$$  Then we have
  \begin{align*}
  \Int_{\Tf^*}(\lf_P|l_{Q})=&\Int_{\Tf^*}^{A}(\lf_P|l_{Q})=2\Int_{\Tf^*}^{A}(\af_{P}^{-}|l_Q)\\
  =&2(\Int_{\Tf^{\bowtie}}^{A}(\af_{P}^{-}|\gamma_{Q}^{+})+\Int_{\Tf^{\bowtie}}^{A}(\af_{P}^{-}|\gamma_{Q}^{-}))\\
  =& \Int_{\Tf^{\bowtie}}^{A}(\af_{P}^{-}|\gamma_{Q}^{-})+\Int_{\Tf^{\bowtie}}^{A}(\af_{P}^{-}|\gamma_{Q}^{+})+\Int_{\Tf^{\bowtie}}^{A}(\af_{P}^{+}|\gamma_{Q}^{-})+\Int_{\Tf^{\bowtie}}^{A}(\af_{P}^{+}|\gamma_{Q}^{+})\\
  =&\Int_{\Tf^{\bowtie}}(\af_{P}^{-}|\gamma_{Q}^{-})+\Int_{\Tf^{\bowtie}}(\af_{P}^{-}|\gamma_{Q}^{+})+\Int_{\Tf^{\bowtie}}(\af_{P}^{+}|\gamma_Q^{+})+\Int_{\Tf^{\bowtie}}(\af_{P}^{+}|\gamma_{Q}^{-}).
   \end{align*}
  where the third equality is due to (\ref{g: a-|l_Q}) ({\em cf.} $\Int(\lf_{P_2}|l_{P_3})$ in Figure \ref{fig:intersection of  l_{p_3} and T}).

\vspace{0.2cm}
\item [Case 2.2]  $P$ is the basepoint of $l_Q$, then $P$ is also the endpoint of both $\gamma_Q^-$ and $\gamma_Q^+$. Since $\gamma_Q^-$ and $\gamma_Q^+$ are compatible and $Q\neq P$, we obtain that $\gamma_Q^-$ and $\gamma_Q^+$ are tagged in the same way at the ends incident to $P$.  Hence, $$\Int_{\Tf^{\bowtie}}^{D}(\af_{P}^{-}|\gamma_{Q}^{-})+\Int_{\Tf^{\bowtie}}^{D}(\af_{P}^{-}|\gamma_{Q}^{+})+\Int_{\Tf^{\bowtie}}^{D}(\af_{P}^{+}|\gamma_{Q}^{-})+\Int_{\Tf^{\bowtie}}^{D}(\af_{P}^{+}|\gamma_{Q}^{+})=2.$$
   Then we have
   \begin{eqnarray*}
      \Int_{\Tf^*}(\lf_P|l_{Q})&=&\Int_{\Tf^*}^{A}(\lf_P|l_{Q})
   =2+2\Int_{\Tf^*}^{A}(\af_{P}^{-}|l_Q)\\
   &=&2(1+\Int_{\Tf^{\bowtie}}^{A}(\af_{P}^{-}|\gamma_{Q}^{+})+\Int_{\Tf^{\bowtie}}^{A}(\af_{P}^{-}|\gamma_{Q}^{+}))\\
   &=&2+\Int_{\Tf^{\bowtie}}^{A}(\af_{P}^{-}|\gamma_{Q}^{-})+\Int_{\Tf^{\bowtie}}^{A}(\af_{P}^{-}|\gamma_{Q}^{+})+\Int_{\Tf^{\bowtie}}^{A}(\af_{P}^{+}|\gamma_{Q}^{-})+\Int_{\Tf^{\bowtie}}^{A}(\af_{P}^{+}|\gamma_{Q}^{+})\\
   &=&\Int_{\Tf^{\bowtie}}(\af_{P}^{-}|\gamma_{Q}^{-})+\Int_{\Tf^{\bowtie}}(\af_{P}^{-}|\gamma_{Q}^{+})+\Int_{\Tf^{\bowtie}}(\af_{P}^{+}|\gamma_{Q}^{-})+\Int_{\Tf^{\bowtie}}(\af_{P}^{+}|\gamma_{Q}^{+}).
   \end{eqnarray*}
 where the third equality is due to (\ref{g:a-|l_p}) ({\em cf.} $\Int(\lf_{P_1}|l_{P_3})$ in Figure \ref{fig:intersection of  l_{p_3} and T}).
\end{itemize}
Therefore, for each  conjugate pair $\{\gamma_{Q}^{-},\gamma_{Q}^{+}\}$ in $\Mcal_1$ with $Q\neq P$, we have 
\begin{equation}\label{g: case wilde{a_p}+ with l_Q}
 \Int_{\Tf^*}(\lf_P|l_{Q})=\Int_{\Tf^{\bowtie}}(\af_{P}^{-}|\gamma_{Q}^{-})+\Int_{\Tf^{\bowtie}}(\af_{P}^{-}|\gamma_{Q}^{+})+\Int_{\Tf^{\bowtie}}(\af_{P}^{+}|\gamma_{Q}^{-})+\Int_{\Tf^{\bowtie}}(\af_{P}^{+}|\gamma_{Q}^{+}). 
    \end{equation}  

\vspace{0.2cm}

For $\{\gamma_{P}^{-},\gamma_{P}^{+}\}$ in $\Mcal_1$,     
we have, $$\Int_{\Tf^*}(\lf_P|l_{P})=\Int^A_{\Tf^*}(\lf_P|l_{P})=2\Int^A_{\Tf^*}(\af_{P}^{-}|l_{P})=2(1+\Int_{\Tf^{\bowtie}}^{A}(\af_{P}^{-}|\gamma^-_{P})+\Int_{\Tf^{\bowtie}}^{A}(\af_{P}^{-}|\gamma^+_{P})),$$
where the third equality is due to (\ref{g:a-|l_p}) ({\em cf.} $\Int(\lf_{P_3}|l_{P_3})$ in Figure \ref{fig:intersection of  l_{p_3} and T}). Moreover, since $$\Int_{\Tf^{\bowtie}}^{D}(\af_{P}^{-}|\gamma_{P}^{-})+\Int_{\Tf^{\bowtie}}^{D}(\af_{P}^{-}|\gamma_{P}^{+})+\Int_{\Tf^{\bowtie}}^{D}(\af_{P}^{+}|\gamma_{P}^{-})+\Int_{\Tf^{\bowtie}}^{D}(\af_{P}^{+}|\gamma_{P}^{+})=2,$$
we obtain 
\begin{equation}\label{g: case wilde{a_p}+ with l_P}
    \Int_{\Tf^*}(\lf_P|l_{P})=\Int_{\Tf^{\bowtie}}^{}(\af_{P}^{-}|\gamma^-_{P})+\Int_{\Tf^{\bowtie}}^{}(\af_{P}^{-}|\gamma^+_{P})+\Int_{\Tf^{\bowtie}}^{}(\af_{P}^{+}|\gamma^-_{P})+\Int_{\Tf^{\bowtie}}^{}(\af_{P}^{+}|\gamma^+_{P}).
\end{equation}
By (\ref{g: case wilde{a_p}+ with gamma_0}),  (\ref{g: case wilde{a_p}+ with l_Q}) and (\ref{g: case wilde{a_p}+ with l_P}), we have $$\Int_{\Tf^*}(\lf_P|\Mcal^*)=\Int_{\Tf^{\bowtie}}(\af_{P}^{-}|\Mcal)+\Int_{\Tf^{\bowtie}}(\af_{P}^{+}|\Mcal).$$

Thus, statement $(3)$ holds.

\vspace{0.2cm}

For (4), because $\af$ has no puncture as endpoint, thus  $\Int_{\Tf^{\bowtie}}^{D}(\af|\gamma)=0$ for any  $\gamma\in\Mcal$. Then  for each   $\gamma\in\Mcal_2$,
\begin{eqnarray*}\label{g: a+ neq a|gamma0}
    \Int_{\Tf^*}(\af|\gamma^{\circ})=\Int_{\Tf^*}^{A}(\af|\gamma^{\circ})=\Int_{\Tf^{\bowtie}}^{A}(\af|\gamma)+\Int_{\Tf^{\bowtie}}^{D}(\af|\gamma)=\Int_{\Tf^{\bowtie}}(\af|\gamma).
\end{eqnarray*}
Moreover, for each  conjugate pair $\{\gamma_{Q}^{-},\gamma_{Q}^{+}\}$ in $\Mcal_1$,
\begin{eqnarray*}\label{g:a|l_Q}
   \Int_{\Tf^*}(\af|l_{Q})
   &=&2\Int_{\Tf^*}^{A}(\af|(\gamma^-_{Q})^\circ)
   =2\Int_{\Tf^{\bowtie}}^{A}(\af|\gamma^-_{Q})\\
     &=&\Int_{\Tf^{\bowtie}}^{A}(\af|\gamma^-_{Q})+\Int_{\Tf^{\bowtie}}^{A}(\af|\gamma^+_{Q}) \\
      &=&\Int_{\Tf^{\bowtie}}^{A}(\af|\gamma^-_{Q})+\Int_{\Tf^{\bowtie}}^{A}(\af|\gamma^+_{Q})+\Int_{\Tf^{\bowtie}}^{D}(\af|\gamma^-_{Q})+\Int_{\Tf^{\bowtie}}^{D}(\af|\gamma^+_{Q}) \\
      &=&\Int_{\Tf^{\bowtie}}(\af|\gamma^-_{Q})+\Int_{\Tf^{\bowtie}}(\af|\gamma^+_{Q}). 
\end{eqnarray*}
Then  we have $\Int_{\Tf^*}(\af|\Mcal^*)= \Int_{\Tf^{\bowtie}}(\af|\Mcal).$ So statement $(4)$ holds.
\end{proof}
Because $\Int_{\Tf^{\bowtie}}^{A}(\af_{P}^{+}|\Mcal)=\Int_{\Tf^{\bowtie}}^{A}(\af_{P}^{-}|\Mcal)$ and $\Int_{\Tf^{\bowtie}}(\af|\Mcal)=\Int_{\Tf^{\bowtie}}^A(\af|\Mcal)+\Int_{\Tf^{\bowtie}}^{D}(\af|\Mcal)$ for any $\af\in\{\af_P^-,\af_P^+\}$, so we have
\begin{eqnarray}\label{g:S_M=plain-notched}
   S_{\Mcal}
   &=&\left\{P\in \Pf\Big| \Int_{\Tf^{\bowtie}}^{D}(\af_{P}^{+}|\Mcal)-\Int_{\Tf^{\bowtie}}^{D}(\af_{P}^{-}|\Mcal)< 0\right\}\\
   &=&\left\{P\in \Pf\Big| \Int_{\Tf^{\bowtie}}^{}(\af_{P}^{+}|\Mcal)-\Int_{\Tf^{\bowtie}}^{}(\af_{P}^{-}|\Mcal)< 0\right\}.\notag
\end{eqnarray}

Now we can state our first main  result.
\begin{theorem}\label{t:equivalence of T^* and T^bowtie}
    Let $(\Sf,\Mf,\Pf,\Tf)$ be a skew-tiling, $\Tf^{\bowtie}$ its tagged version, and   $(\Sf^*,\Mf^*,\Tf^*)$ the unfolded tiling of $(\Sf,\Mf,\Pf,\Tf)$. Then for any $\Mcal,\Ncal\in\mathscr{R}(\Sf)$, the following are equivalent:
\begin{itemize}
    \item [(1)] $\Intv_{\Tf^{\bowtie}}(\Mcal)=\Intv_{\Tf^{\bowtie}}(\Ncal)$;\vspace{0.15cm}
    \item [(2)]  $\Intv_{\Tf^*}(\Mcal^*)=\Intv_{\Tf^*}(\Ncal^*)$  and  $S_{\Mcal}=S_{\Ncal}.$
\end{itemize}
    \end{theorem}
\begin{proof}
$(1)\Rightarrow (2)$. 
  Since $\Intv_{\Tf^{\bowtie}}(\Mcal) = \Intv_{\Tf^{\bowtie}}(\Ncal)$, 
we have $\Int_{\Tf^{\bowtie}}^{}(\af|\Mcal)=\Int_{\Tf^{\bowtie}}^{}(\af|\Ncal)$ for any $\af\in \Tf^{\bowtie}$. In particular, for each puncture $P\in\Pf$, 
\begin{eqnarray*}
    \Int_{\Tf^{\bowtie}}^{}(\af_{P}^{+}|\Mcal)-\Int_{\Tf^{\bowtie}}^{}(\af_{P}^{-}|\Mcal)
    =\Int_{\Tf^{\bowtie}}^{}(\af_{P}^{+}|\Ncal)-\Int_{\Tf^{\bowtie}}^{}(\af_{P}^{-}|\Ncal).
\end{eqnarray*}
 By (\ref{g:S_M=plain-notched}), 
 $S_{\Mcal}=S_{\Ncal}$.
 Let $\af\in\Tf^*$. We consider the following cases:
\begin{itemize}
    \item[Case 1.1:] $\af=\lf_P$ for some puncture $P$.  By Lemma \ref{l:skew-gentle to gentle}(3), 
\begin{eqnarray*}
    \Int_{\Tf^*}(\lf_P|\Mcal^*)
   =\Int_{\Tf^{\bowtie}}(\af_{P}^{-}|\Mcal)+\Int_{\Tf^{\bowtie}}(\af_{P}^{+}|\Mcal)
   &=&\Int_{\Tf^{\bowtie}}(\af_{P}^{-}|\Ncal)+\Int_{\Tf^{\bowtie}}(\af_{P}^{+}|\Ncal)\\
   &=&\Int_{\Tf^*}(\lf_P|\Ncal^*). 
\end{eqnarray*} 
   \item[Case 1.2:] $\af=\af_P^{-}$ for some puncture $P$. If $P\in S_{\Mcal}$, then $P\in S_{\Ncal}$, and by Lemma \ref{l:skew-gentle to gentle}(2),  $$\Int_{\Tf^*}(\af_{P}^{-}|\Mcal^*)=\Int_{\Tf^{\bowtie}}(\af_{P}^{-}|\Mcal)=\Int_{\Tf^{\bowtie}}(\af_{P}^{-}|\Ncal)=\Int_{\Tf^*}(\af_{P}^{-}|\Ncal^*).$$ Otherwise, if $P\notin S_{\Mcal}$, then $P\notin S_{\Ncal}$, and by Lemma \ref{l:skew-gentle to gentle}(1), 
    $$\Int_{\Tf^*}(\af_{P}^{-}|\Mcal^*)=\Int_{\Tf^{\bowtie}}(\af_{P}^{+}|\Mcal)=\Int_{\Tf^{\bowtie}}(\af_{P}^{+}|\Ncal)=\Int_{\Tf^*}(\af_{P}^{-}|\Ncal^*).$$
    \item[Case 1.3:] $\af\neq \lf_P$ for any puncture $P$ and $\af\neq\af_Q^{-}$ for any puncture $Q$.  By Lemma \ref{l:skew-gentle to gentle}(4),
    $$\Int_{\Tf^*}(\af|\Mcal^*)=\Int_{\Tf^{\bowtie}}(\af|\Mcal)=\Int_{\Tf^{\bowtie}}(\af|\Ncal)=\Int_{\Tf^*}(\af|\Ncal^*).$$
\end{itemize}
Hence, we conclude that 
$\Intv_{\Tf^*}(\Mcal^*)=\Intv_{\Tf^*}(\Ncal^*).$

$(2)\Rightarrow (1)$. For any puncture $P$,  by Lemma \ref{l:skew-gentle to gentle}(3), we have $$ \Int_{\Tf^*}(\lf_P|\Mcal^*)=
    \Int_{\Tf^{\bowtie}}(\af_{P}^{-}|\Mcal)+\Int_{\Tf^{\bowtie}}(\af_{P}^{+}|\Mcal),$$
    and $$ \Int_{\Tf^*}(\lf_P|\Ncal^*)=
    \Int_{\Tf^{\bowtie}}(\af_{P}^{-}|\Ncal)+\Int_{\Tf^{\bowtie}}(\af_{P}^{+}|\Ncal).$$
 Since  $\Intv_{\Tf^*}(\Mcal^*)=\Intv_{\Tf^*}(\Ncal^*)$, it follows that 
 \begin{equation}\label{g:m_-+m_+=n_-+n_+}
      \Int_{\Tf^{\bowtie}}(\af_{P}^{-}|\Mcal)+\Int_{\Tf^{\bowtie}}(\af_{P}^{+}|\Mcal) =\Int_{\Tf^{\bowtie}}(\af_{P}^{-}|\Ncal)+\Int_{\Tf^{\bowtie}}(\af_{P}^{+}|\Ncal).
 \end{equation}
Now let $\af\in\Tf^{\bowtie}$, we consider the following cases:
\begin{itemize}
    \item[Case 2.1:] $\af\in\{\af_{P}^{-},\af_{P}^{+}\}$ for some puncture $P\notin S_{\Mcal}$.
Then $P\notin S_{\Ncal}$. By Lemma \ref{l:skew-gentle to gentle}$(1)$,
 \begin{equation*}
  \Int_{\Tf^{\bowtie}}(\af_{P}^{-}|\Mcal)=\Int_{\Tf^*}(\af_{P}^{-}|\Mcal^*)={\Int_{\Tf^*}(\af_{P}^{-}|\Ncal^*)=\Int_{\Tf^{\bowtie}}(\af_{P}^{-}|\Ncal)}.
  \end{equation*}
  Combining with (\ref{g:m_-+m_+=n_-+n_+}) yields
$\Int_{\Tf^{\bowtie}}(\af_{P}^{+}|\Mcal)=\Int_{\Tf^{\bowtie}}(\af_{P}^{+}|\Ncal).$
  \item[Case 2.2:] $\af\in\{\af_{P}^{-},\af_{P}^{+}\}$  for some puncture $P\in S_{\Mcal}$. Then $P\in S_\Ncal$. By Lemma \ref{l:skew-gentle to gentle}$(2)$, 
\begin{equation*}\Int_{\Tf^{\bowtie}}(\af_{P}^{+}|\Mcal)=\Int_{\Tf^*}(\af_{P}^{-}|\Mcal^*)=\Int_{\Tf^*}(\af_{P}^{-}|\Ncal^*)=\Int_{\Tf^{\bowtie}}(\af_{P}^{+}|\Ncal).
\end{equation*}
 Combining with (\ref{g:m_-+m_+=n_-+n_+}) gives
$\Int_{\Tf^{\bowtie}}(\af_{P}^{-}|\Mcal)=\Int_{\Tf^{\bowtie}}(\af_{P}^{-}|\Ncal).$ 
\item[Case 2.3:]  No puncture is an endpoint of $\af$.  By Lemma \ref{l:skew-gentle to gentle}(4), $$\Int_{\Tf^{\bowtie}}(\af|\Mcal)=\Int_{\Tf^*}(\af|\Mcal^*)=\Int_{\Tf^*}(\af|\Ncal^*)=\Int_{\Tf^{\bowtie}}(\af|\Ncal).$$
\end{itemize}
\noindent Hence, we conclude that $\Intv_{\Tf^{\bowtie}}(\Mcal)=\Intv_{\Tf^{\bowtie}}(\Ncal).$
\end{proof}

     
\subsection{Multisets determined by their intersection vectors }\label{s:subs3}
\begin{theorem}\cite[Theorem 1.1]{FG22}\label{t:no even cycle}
Let $(\Sf,\Mf,\Tf)$ be a tiling. Suppose $(\Sf,\Mf,\Tf)$ contains no tile of type $\II$ nor even-gons of type $\V$. Then for any $\Mcal,\ \Ncal\in\mathscr{R}(\Sf)$ satisfying $\underline{\Int}_{\Tf}(\Mcal)=\underline{\Int}_{\Tf}(\Ncal)$, we have $\Mcal=\Ncal$.
\end{theorem}

Combining Theorem \ref{t:equivalence of T^* and T^bowtie} and Theorem \ref{t:no even cycle}, we have
\begin{theorem}\label{t:tagged multiset determines}
Let $(\Sf,\Mf,\Pf,\Tf)$ be a skew-tiling that   contains no tile of type $\II$ or even gons of type $\V$. Let $\Tf^{\bowtie}$ be the tagged version of $\Tf$. Then for any $\Mcal,\ \Ncal\in\mathscr{R}(\Sf)$ satisfying $\underline{\Int}_{\Tf^{\bowtie}}(\Mcal)=\underline{\Int}_{\Tf^{\bowtie}}(\Ncal)$, we have $\Mcal=\Ncal$. 
\end{theorem}
\begin{proof}
Let $(\Sf^*,\Mf^*,\Tf^*)$ be the associated unfolded tiling of  $(\Sf,\Mf,\Pf,\Tf)$.
        By Theorem \ref{t:equivalence of T^* and T^bowtie}, we have $\Int_{\Tf^*}(\Mcal^*)=\Int_{\Tf^*}(\Ncal^*)$ and $S_\Mcal=S_\Ncal$. Since in the skew-tiling  $(\Sf,\Mf,\Pf,\Tf)$, there are no tiles of type $\II$ or even-gons of type $\V$. Thus, by the construction of  $(\Sf^*,\Mf^*,\Tf^*)$, it is easy to get that $(\Sf^*,\Mf^*,\Tf^*)$ also has no tiles of type $\II$ or even-gons of type $\V$. Therefore, by applying Theorem \ref{t:no even cycle}, we conclude that $\Mcal^* = \Ncal^*$.
        
By Subsection~\ref{ss:Permissible arcs from  unfolded tiling to skew-tiling}, we know that $\Mcal$ (resp. $\Ncal$) is uniquely determined by $\Mcal^*$ and $S_{\Mcal}$ (resp. $S_\Ncal$). Thus we obtain $\Mcal=\Ncal$.
\end{proof}

\section{Application to skew-gentle algebras}\label{s:s4}


\subsection{Tagged permissible arcs and $\tau$-rigid modules} Let $A$ be a finite dimensional algebra, $\{e_1,\cdots,e_n\}$  a complete set of primitive orthogonal idempotents of $A$. For a right $A$-module $M$, denote by $\dimv M=(\dim_K Me_1,\cdots,\dim_KMe_n)^{tr}$ the dimension vector of $M$, where $(-)^{tr}$ denotes the transpose of $(-)$.
Let $\tau$ be the Auslander-Reiten translation. 
A finitely generated $A$-module $M$ is called {\it $\tau$-rigid} if $\Hom_A(M,\tau M)=0$. Clearly, every finitely generated projective $A$-module is $\tau$-rigid. 
Let ${\tau}\opname{-rigid} A$ denote the set of isomorphism classes of $\tau$-rigid modules in $\mathrm{mod} A$, and let $\opname{ind}\tau\opname{-rigid} A$ denote the subset of indecomposable $\tau$-rigid $A$-modules. 
The following results are consequences of those in \cite{HZZ}.
\begin{lemma}\cite[Theorems 2.29, 3.12 and 4.3, Propositions 5.3 and 5.6]{HZZ}\label{l:intersection and dimension}
Let $A=K Q/\langle I\rangle$ be a finite-dimensional skew-gentle algebra. Then there exists a skew-tiling $(\Sf,\Mf,\Pf,\Tf)$ with the corresponding admissible partial tagged triangulation $\Tf^{\bowtie}$ such that 
there is a bijection
			\begin{eqnarray*}
			M\colon &\{\gamma\mid \gamma\in\PTAS\}&\to\opname{ ind~\tau-}\opname{rigid} A,\\
			&\gamma&\mapsto M(\gamma)
			\end{eqnarray*}
			satisfying $\underline{\Int}_{\Tf^{\bowtie}}(\gamma)=\dimv M(\gamma)$ for any  $\gamma\in\PTAS$.
	
  Moreover, 
  the bijection $M$ induces a bijection
		\[\mathscr{R}(\Sf)\to\opname{\tau-}\opname{rigid} A,\]
		mapping $\mathscr{R}\in \mathscr{R}(\Sf)$ to $\bigoplus\limits_{\gamma\in \mathscr{R}}M(\gamma)$  satisfying $\underline{\Int}_{\Tf^{\bowtie}}(\mathscr{R})=\dimv \bigoplus\limits_{\gamma\in \mathscr{R}}M(\gamma)$.
\end{lemma}


\subsection{Cartan matrices between skew-tiling and unfolded tiling}
For a finite dimensional algebra $A=KQ/\langle I\rangle$, an oriented cycle $c=\alpha_1\dots\alpha_s$ of length $s$ in $Q$ is said to have \emph {full zero relations} if $\alpha_i\alpha_{i+1}\in I$ for all $i=1,\dots s-1$ and $\alpha_s\alpha_1 \in I.$ Such a cycle is called \emph {minimal} if the arrows $\alpha_1,\dots,\alpha_s$ are pairwise different. Denote  by $ec(Q,I)$ the number of  minimal oriented cycles of even length in $Q$ having full zero relations, and by $oc(Q,I)$ the number of  minimal oriented cycles of odd length in $Q$ having full zero relations.
 
 Let $\{e_1,\cdots,e_n\}$ be a complete set of primitive  orthogonal idempotents of $A$. For each $1\leq i\leq n$, let $P_i=e_iA$ be the corresponding projective right $A$-module. The {\it Cartan matrix} $\mathbf{C}_A$ of $A$ is defined as 
 \[\mathbf{C}_A=(\dimv P_1,\dots,\dimv P_n).
\] The determinant $\det \mathbf{C}_A$ of $\mathbf{C}_A$ has been explicitly computed for both gentle and skew-gentle algebras.
\begin{lemma}$\mathrm{(}$\cite[Theorem 1]{H} and \cite[Theorem 4.1]{BH}$\mathrm{)}$\label{l:det}
 Let $(Q,Sp,I)$ be a skew-gentle triple. 
 $A^g = KQ/\langle I\rangle$ the associated gentle algebra, $A^{\bowtie}= KQ^{\bowtie}/ \langle I^{\bowtie}\rangle$ its associated skewed-gentle algebra.  Then 
 \[\det \mathbf{C}_{A^{\bowtie}}=\det \mathbf{C}_{A^g}=\begin{cases}
     0  & \text{ if } ec(Q,I)>0;\\
     2^{oc(Q,I)} &\text{ otherwise}.
 \end{cases}\]
\end{lemma}

Now we have
\begin{proposition}\label{p:det between skewed-gentle and unfolded tiling}
Let  $(\Sf,\Mf,\Pf)$ be a marked surface with an skew-tiling $\Tf$, $\Tf^{\bowtie}$ the tagged version of $\Tf$, $(\Sf^*,\Mf^*,\Tf^*)$ the associated {unfolded tiling} and $(Q,Sp,I)$  the associated skew-gentle triple. Then $$\det \mathbf{C}_{A^{{\bowtie}}}=\det \mathbf{C}_{A_{\Tf^{*}}}.$$
\end{proposition}
\begin{proof}
By Lemma \ref{l:det}, we have  $\det \mathbf{C}_{A^{{\bowtie}}}=\det \mathbf{C}_{A^g}$.  
 According to Proposition~\ref{p:relation between tiling algebras}, $A_{\Tf^*}$ is isomorphic to $A^*=kQ^*/\langle I^*\rangle$, where
\begin{itemize}
    \item $Q^*_0=Q_0\cup \{i^{*}\mid i\in Sp\}$,  \item $Q^*_1=Q_1\cup\{\rho_i:i\rightarrow i^*, \rho_{i^*}:i^*\rightarrow i\mid i\in Sp\}$;
     \vspace{0.1cm}
    \item $I^*=I\cup \{\rho_{i}\rho_{i^*}\mid i\in Sp\}$.
\end{itemize} 
Then one can get that $oc(Q,I)=oc(Q^*,I^*)$ and $ec(Q,I)=ec(Q^*,I^*)$.
Thus, by Lemma \ref{l:det}, we have $$\det \mathbf{C}_{A^{{\bowtie}}}=\det \mathbf{C}_{A^{g}}=\det \mathbf{C}_{A_{\Tf^{*}}}.$$
\end{proof}

\subsection{$\tau$-rigid modules over skew-gentle algebras}
Now we have the following equivalent statements about when the $\tau$-rigid modules over a skew-gentle algebra are uniquely determined by their dimension vectors.
\begin{theorem}\label{t: ns-conditions of tau-rigid}
  Let $(Q,Sp,I)$ be a skew-gentle triple, and   
  let
\[
A^{sg}=KQ^{sp}/\langle I^{sg}\rangle,\quad
A^{\bowtie}=KQ^{\bowtie}/\langle I^{\bowtie}\rangle,\quad
A^{*}=KQ^{*}/\langle I^{*}\rangle
\]
be the associated skew-gentle, skewed-gentle and unfolded tiling algebras, respectively.
The following statements are equivalent:
\begin{itemize}
\item[(1)] Different $\tau$-rigid $A^{sg}$-modules have different dimension vectors;
    \item[$(1')$] Different $\tau$-rigid $A^{\bowtie}$-modules have different dimension vectors;
\item [(2)]The determinant of the Cartan matrix of $A^{sg}$ is nonzero;
\item [$(2')$]The determinant of the Cartan matrix of $A^{\bowtie}$ is nonzero;
\item[(3)] $(Q,I)$
contains no minimal oriented cycle of even length with full zero relations;
\item[(4)]$(Q^{sp},I^{sg})$
contains no minimal oriented cycle of even length with full zero relations;

\item[(5)] Different $\tau$-rigid $A^*$-modules have different dimension vectors.
\end{itemize}

\end{theorem}
\begin{proof}
By Theorem~\ref{t:skew-tiling algebra}, there is a skew-tiling
$(\Sf, \Mf,\Pf,\Tf)$ such that $A^{sg}\cong A_{\Tf}$. Let $\Tf^{\bowtie}$ be the tagged version of $\Tf$ and  $(\Sf^*,\Mf^*,\Tf^*)$ be the associated unfolded tiling. Thus by Proposition~\ref{p:relation between tiling algebras}, we know   $A^*=A_{T^*}$.

``$(1)\Leftrightarrow (1')$"  and ``$(2)\Leftrightarrow (2')$" are  due to the Morita equivalence between  $A^{sg}$ and $A^{\bowtie}$.

``$(1)\Rightarrow (2)$" is clear since projective modules are $\tau$-rigid.

 ``$(2')\Leftrightarrow (3)$" is obtained from Lemma~\ref{l:det}.

  ``$(3)\Leftrightarrow(4)$" follows directly from the construction of $(Q^{sp}, I^{sg})$ and $(Q,I)$.

For ``$(4)\Rightarrow (1)$": let $X$ and $Y$ be two $\tau$-rigid $A^{sg}$-modules such that $\dimv X=\dimv Y$. Then by Lemma \ref{l:intersection and dimension}, there are  two multisets
 $\Mcal,~\Ncal\in\mathscr{R}(\Sf)$ such that $\Intv_{\Tf^{\bowtie}}(\Mcal)=\dimv X$ and $\Intv_{\Tf^{\bowtie}}(\Ncal)=\dimv Y,$ which means
 $\Intv_{\Tf^{\bowtie}}(\Mcal)=\Intv_{\Tf^{\bowtie}}(\Ncal).$  Because $(Q^{sp}, I^{sg})$
contains no minimal oriented cycle of even length with full zero relations, thus $(\Sf, \Mf,\Pf,\Tf)$  has no tiles of type (II) or even-gons of type (V). By Theorem~\ref{t:tagged multiset determines}, $\Mcal=\Ncal$, which means $X=Y$.

 For ``$(4)\Leftrightarrow (5)$":  by the construction of $(\Sf^*,\Mf^*,\Tf^*)$, it is clear that  $(\Sf, \Mf,\Pf,\Tf)$  has no tiles of type (II) or even-gons of type (V) if and only if $(\Sf^*,\Mf^*,\Tf^*)$
 has no tiles of type (II) or even-gons of type (V). On the other hand, $(\Sf^*,\Mf^*,\Tf^*)$
 has no tiles of type (II) or even-gons of type (V) if and only if $(Q^*,I^*)$ contains no minimal oriented cycle of even length with full zero relations, which is equivalent to different $\tau$-rigid $A^*$-modules have different dimension vectors by  \cite[Theorem 1.2]{FG22}.
\end{proof}
If $A$ has finite global dimension, it is clear that $\det \mathbf{C}_A\neq 0$. Then by Theorem~\ref{t: ns-conditions of tau-rigid}, we have
\begin{corollary}
     Let $A$ be a skew-gentle algebra with finite global dimension.
Then different $\tau$-rigid $A$-modules have different dimension vectors.
\end{corollary}

\section{Two examples}\label{s:two examples}
\begin{example}\label{ex:case for odd}
Let $(Q,Sp,I)$ be a skew-gentle triple, where

$\begin{tikzcd}[sep=tiny]
	&&& 1 \\
	{Q=}&&&&&,&&&{I=\{\beta\gamma,\gamma\alpha,\alpha\beta\}},&{Sp=\{1\}.} \\
	& 2 &&&& 3
	\arrow["\alpha"', from=1-4, to=3-2]
	\arrow["\gamma"', from=3-2, to=3-6]
	\arrow["\beta"', from=3-6, to=1-4]
\end{tikzcd}$

 Then

\vspace{0.3cm}

       $\begin{tikzcd}[sep=tiny]
	&&& 1 \\
	{Q^{sp}=} &&&&&,&{I^{sp}=\{\beta\gamma,\gamma\alpha,\alpha\beta,\epsilon_1^2\},}&{I^{sg}=\{\beta\gamma,\gamma\alpha,\alpha\beta,\epsilon_1-\epsilon_1^2\},}\\
	& 2 &&&& 3
	\arrow["{\epsilon_1}", from=1-4, to=1-4, loop, in=55, out=125, distance=10mm]
	\arrow["\alpha"', from=1-4, to=3-2]
	\arrow["\gamma"', from=3-2, to=3-6]
	\arrow["\beta"', from=3-6, to=1-4]
\end{tikzcd}$

\vspace{0.3cm}

 $\begin{tikzcd}[sep=tiny]
	&&& {1^+} \\
	\\
	{Q^{\bowtie}=} & 2 &&&& {3,} && {I^{\bowtie}=\{\alpha^+\beta^+-\alpha^-\beta^-,~\beta^{\pm}\gamma,~\gamma\alpha^{\pm}\},} \\
	\\
	&&& {1^-}
	\arrow["{{{{\alpha^+}}}}"', from=1-4, to=3-2]
	\arrow["\gamma"', from=3-2, to=3-6]
	\arrow["{{{{\beta^+}}}}"', from=3-6, to=1-4]
	\arrow["{{{{\beta^-}}}}", from=3-6, to=5-4]
	\arrow["{{{{\alpha^-}}}}", from=5-4, to=3-2]
\end{tikzcd}$

\vspace{0.3cm}

$\begin{tikzcd}[sep=tiny]
	&&& {1^*} \\
	\\
	&&& 1 \\
	{Q^*=} &&&&& {,}&& & {I^*=\{\beta\gamma,\gamma\alpha,\alpha\beta,\rho_1\rho_{1^*}\},} \\
	& 2 &&&& 3
	\arrow["{{{\rho_{1^*}}}}", shift left, from=1-4, to=3-4]
	\arrow["{{{\rho_1}}}", shift left, from=3-4, to=1-4]
	\arrow["\alpha"', from=3-4, to=5-2]
	\arrow["\gamma"', from=5-2, to=5-6]
	\arrow["\beta"', from=5-6, to=3-4]
\end{tikzcd}$

and 

\vspace{0.3cm}

    $
\mathbf{C}_{A^g}=\begin{bmatrix}
1 & 0 & 1 \\
1 & 1 & 0 \\
0 & 1 & 1
\end{bmatrix},
\qquad 
\mathbf{C}_{A^{sp}}=\begin{bmatrix}
2 & 0 & 2 \\
2 & 1 & 1 \\
0 & 1 & 1
\end{bmatrix},
$

\vspace{0.3cm}

$
\mathbf{C}_{A^{sg}}=\mathbf{C}_{A^{\bowtie}}=(\dimv P_{1^-},\dimv P_{1^{+}},\dimv P_2,\dimv P_3)=\begin{bmatrix}
1 & 0 & 0 & 1 \\
0 & 1 & 0 & 1 \\
1 & 1 & 1 & 1 \\
0 & 0 & 1 & 1 \\
\end{bmatrix},$

\vspace{0.3cm}

$
\mathbf{C}_{A^*}=(\dimv P_1,\dimv P_{1^*},\dimv P_2,\dimv P_3)=\begin{bmatrix}
2 & 1 & 0 & 2 \\
1 & 1 & 0 & 1 \\
2 & 1 & 1 & 1 \\
0 & 0 & 1 & 1 \\
\end{bmatrix}.
$

\vspace{0.3cm}

\noindent where  $A^{g}=KQ^{}/\langle I^{}\rangle$, $A^{sp}=KQ^{sp}/\langle I^{sp}\rangle$ and $A^*=KQ^*/\langle I^*\rangle$ are all gentle algebras, 
$A^{sg}=KQ^{sp}/\langle I^{sg}\rangle$ is the associated skew-gentle algebra, 
$A^{\bowtie}=KQ^{\bowtie}/\langle I^{\bowtie}\rangle$ is the associated skewed-gentle algebra. In particular, $A^*=KQ^*/\langle I^*\rangle$ is the associated unfolded tiling  algebra.
It is easy to get that  $$\det \mathbf{C}_{A^g}=\det\mathbf{C}_{A^{\bowtie}}=\det\mathbf{C}_{A^*}=2,\text{~and~}\det\mathbf{C}_{A^{sp}}=4.$$
By Theorem~\ref{t: ns-conditions of tau-rigid},  the $\tau$-rigid  $A^{sg}$-modules and the $\tau$-rigid  $A^*$-modules are uniquely determined by their dimension vectors.

\end{example}
\begin{example}\label{ex:case for even}
Compared to Example \ref{ex:case for odd}, now we consider another skew-gentle triple  $(Q,Sp,I)$, where 

$
\begin{tikzcd}[sep=tiny,baseline=(current bounding box.west)]
	& 1 &&& 2 \\
	\\
	{Q=} &&&&& {,} && {I=\{\alpha_2\alpha_1,\alpha_3\alpha_2,\alpha_4\alpha_3,\alpha_1\alpha_4\},} & {Sp=\{1\}}. \\
	& 4 &&& 3
	\arrow["{\alpha_1}", from=1-2, to=1-5]
	\arrow["{\alpha_2}", from=1-5, to=4-5]
	\arrow["{\alpha_4}", from=4-2, to=1-2]
	\arrow["{\alpha_3}", from=4-5, to=4-2]
\end{tikzcd}
$

Then\begin{flushleft}
\begin{tikzcd}[sep=tiny,baseline=(current bounding box.west)]
	& 1 &&& 2 \\
	\\
	{Q^{sp}=} &&&&& {,} && {I^{sp}=\{\alpha_2\alpha_1,\alpha_3\alpha_2,\alpha_4\alpha_3,\alpha_1\alpha_4,\epsilon_1^2\},} &  \\
	& 4 &&& 3 &&& {I^{sg}=\{\alpha_2\alpha_1,\alpha_3\alpha_2,\alpha_4\alpha_3,\alpha_1\alpha_4,\epsilon_1-\epsilon_1^2\},}
	\arrow["{{{\epsilon_1}}}"{pos=0.2}, from=1-2, to=1-2, loop, in=55, out=125, distance=10mm]
	\arrow["{{{\alpha_1}}}", from=1-2, to=1-5]
	\arrow["{{{\alpha_2}}}", from=1-5, to=4-5]
	\arrow["{{{\alpha_4}}}", from=4-2, to=1-2]
	\arrow["{{{\alpha_3}}}", from=4-5, to=4-2]
\end{tikzcd}\end{flushleft}
\begin{flushleft}
\begin{tikzcd}[sep=tiny,baseline=(current bounding box.west)]
	&&& {1^+} &&& 2 \\
	&&&& {1^-} \\
	{Q^{\bowtie}=} &&&&&&& {,} && {I^{\bowtie}=\{\alpha_2\alpha_1^{\pm1},\alpha_3\alpha_2,\alpha_4^{\pm1}\alpha_3,\alpha_1^{+}\alpha_4^+-\alpha_1^{-}\alpha_4^{-}\},} \\
	&&& 4 &&& 3
	\arrow["{\alpha_1^+}", from=1-4, to=1-7]
	\arrow["{\alpha_2}", from=1-7, to=4-7]
	\arrow["{\alpha_1^-}"', from=2-5, to=1-7]
	\arrow["{\alpha_4^+}", from=4-4, to=1-4]
	\arrow["{\alpha_4^-}"', from=4-4, to=2-5]
	\arrow["{\alpha_3}", from=4-7, to=4-4]
\end{tikzcd}\end{flushleft}
\begin{flushleft}
\begin{tikzcd}[sep=tiny,baseline=(current bounding box.west)]
	& {1^*} && 1 &&& 2 \\
	\\
	{Q^*=} &&&&&&& {,} && {I^*=\{\alpha_2\alpha_1,\alpha_3\alpha_2,\alpha_4\alpha_3,\alpha_1\alpha_4,\rho_1\rho_{1^*}\},} \\
	&&& 4 &&& 3
	\arrow["{\rho_{1^*}}", shift left=2, from=1-2, to=1-4]
	\arrow["{\rho_{1}}", shift left=2, from=1-4, to=1-2]
	\arrow["{\alpha_1}", from=1-4, to=1-7]
	\arrow["{\alpha_2}", from=1-7, to=4-7]
	\arrow["{\alpha_4}", from=4-4, to=1-4]
	\arrow["{\alpha_3}", from=4-7, to=4-4]
\end{tikzcd}\end{flushleft}
and

$
\mathbf{C}_{A^g}=\begin{bmatrix}
1 & 0 & 0 & 1\\
1 & 1 & 0 & 0\\
0 & 1 & 1 & 0\\
0 & 0 & 1 & 1\\
\end{bmatrix},
\quad \mathbf{C}_{A^{sp}}=\begin{bmatrix}
2 & 0 & 0 & 2\\
2 & 1 & 0 & 1\\
0 & 1 & 1 & 0\\
0 & 0 & 1 & 1\\
\end{bmatrix}, $

\vspace{0.3cm}

$
\mathbf{C}_{A^{sg}}=\mathbf{C}_{A^{\bowtie}}=(\dimv P_{1^-},\dimv P_{1^{+}},\dimv P_2,\dimv P_3, \dimv P_4)=\begin{bmatrix}
1 & 0 & 0 & 0 &1\\
0 & 1 & 0 & 0 &1 \\
1 & 1 & 1 & 0 &1\\
0 & 0 & 1 & 1 &0\\
0 & 0 & 0 & 1 &1\\
\end{bmatrix},
$

\vspace{0.3cm}

$
 \mathbf{C}_{A^{*}}=(\dimv P_1,\dimv P_{1^*},\dimv P_2,\dimv P_3,\dimv P_4)=\begin{bmatrix}
2 & 1 & 0 & 0 &2\\
1 & 1 & 0 & 0 &1\\
2 & 1 & 1 & 0 &1\\
0 & 0 & 1 & 1 &0\\
0 & 0 & 0 & 1 &1\\
\end{bmatrix}.  
$

\vspace{0.3cm}

\noindent It is easy to get $$\det\mathbf{C}_{A^{g}}=\det\mathbf{C}_{A^{sp}}=\det\mathbf{C}_{A^{\bowtie}}=\det\mathbf{C}_{A^*}=0,$$ 
Moreover, one can get that for $A^{sg}$, $$\dimv(P_{1^-}\oplus P_{1^+}\oplus P_3)=\dimv(P_2\oplus P_4).$$
For $A^*$, $$\dimv(P_1\oplus P_3)=\dimv(P_2\oplus P_4).$$
\end{example}

    \bibliographystyle{acm}
    \bibliography{skewref}
    
\end{document}